\documentclass[12pt]{article}

\advance \topmargin by -\headheight
\advance \topmargin by -\headsep

\evensidemargin \oddsidemargin
\usepackage[dvips]{graphics,color}
\usepackage {pstcol}

\newpsobject{showgrid}{psgrid}{subgriddiv=1,griddots=10,gridlabels=6pt}
\usepackage{amssymb}
\newcommand{\nc}{\newcommand}
\def\pf{\textrm{pf}\,}

\def\om{\omega}

\nc\intint{\int\!\!\int}
\nc{\twotwo}[4]{\left(\begin{array}{cc}#1&#2\\&\\#3&#4\end{array}\right)}
\nc{\twoone}[2]{\left(\begin{array}{c}#1\\\\#2\end{array}\right)}
 \nc{\wh}{\widehat}
 \nc{\pl}{\partial}
 \renewcommand{\sp}{\vskip1ex}
 \nc{\inv}{^{-1}}
\def\ot{\otimes}

\def\La{\Lambda}

\def\ra{\rightarrow}

\def\iy{\infty}
\newcommand{\mi}{\mathrm{i}}
\newcommand{\me}{\mathrm{e}}
\def\hf{{1\over 2}}

\def\be{\begin{equation}}
\def\ee{\end{equation}}
\def\ba{\begin{eqnarray*}}
\def\ea{\end{eqnarray*}}
\def\bae{\begin{eqnarray}}
\def\eae{\end{eqnarray}}
\def\bc{\begin{center}}
\def\ec{\end{center}}
\def\ov{\over}
\def\df{\mathrm{d}}
\def\al{\alpha}
\def\ga{\gamma}
\def\Ga{\Gamma}
\def\s{\sigma}
\def\om{\omega}
 \renewcommand{\th}{\theta}
\newcommand{\And}{ \ \ \textrm{and}\ \  }
\def\z{\zeta}
\def\ve{\varepsilon}
\def\vp{\varphi}

\def\tQ{\widetilde Q}

\def\la{\lambda}
\def\pr{\textrm{P}}

\def\dl{\delta}

\def\x{\xi}
\def\cP{\mathcal{P}}
\def\cD{\mathcal{D}}
\def\cS{\mathcal{S}}

\def\cK{\mathcal{K}}

\def\bN{\mathbb{N}}
\def\bP{\mathbb{P}}

\def\bQ{\mathbb{Q}}
\def\bZ{\mathbb{Z}}
\def\1p{1^\prime}
\def\2p{2^\prime}
\def\3p{3^\prime}
\def\4p{4^\prime}
\def\5p{5^\prime}
\def\6p{6^\prime}
\def\7p{7^\prime}
\def\8p{8^\prime}
\def\9p{9^\prime}
\begin{document}
\title{A Limit Theorem for Shifted Schur Measures}
\author{Craig A.~Tracy\\
Department of Mathematics \\
University of California\\
Davis, CA 95616\\
email: tracy@math.ucdavis.edu\\
\quad\\
 and \\
\quad\\
Harold Widom\\
Department of Mathematics\\
University of California\\
Santa Cruz, CA 95064\\
email: widom@math.ucsc.edu}
\maketitle
\begin{abstract}To each partition $\la=(\la_1,\la_2,\ldots)$  with distinct parts we assign the
probability $Q_\la(x) P_\la(y)/Z$
where $Q_\la$ and $P_\la$ are the Schur $Q$-functions and $Z$ is a normalization
constant.  This measure, which we call the shifted Schur measure, is  analogous
to the much-studied Schur measure.  For the specialization of
the first $m$ coordinates of $x$ and the first $n$ coordinates of $y$ equal to
$\al$ ($0<\al<1$) and the rest equal to zero, we derive a limit law for $\la_1$ as
$m,n\ra\iy$ with $\tau=m/n$ fixed.  For the Schur measure the $\al$-specialization limit
law was derived by Johansson.  Our main
result implies that the two limit laws are identical.
\end{abstract}
\newtheorem{theorem}{Theorem}
\newpage
\setcounter{tocdepth}{2}
\tableofcontents
\section{Introduction}
The Schur measure~\cite{Ok1} assigns to each partition $\la=(\la_1,\la_2,\ldots )$ the
weight
\[  s_\la(x) \, s_\la(y) \]
where $s_\la$ are the Schur functions. (See, e.g.,~\cite{Mac, St}.)  Thus
\be \sum_{{\la\in\cP\atop \la_1\le h}} s_\la(x) \, s_\la(y)  \label{GesselLHS}\ee
is the (unnormalized) probability that $\la_1$, the number of boxes in the first row of the
associated Young diagram, is less than or equal to $h$.  Here $\cP$ denotes the set
of all partitions. The normalization constant, $Z$, is determined from the
$h\ra\iy$ limit
\be Z:=\sum_{\la\in\cP} s_\la(x)\, s_\la(y)=\prod_{i, j}{1\ov 1-x_i y_j},\, \label{Cauchy}\ee
where the last equality is the Cauchy identity for Schur functions.
A theorem of Gessel~\cite{Ge}
expresses the partition sum (\ref{GesselLHS})  as an $h\times h$ Toeplitz determinant $D_h(\vp)$.
It follows from this that the normalization constant is also given by
\[ Z=\lim_{h\ra\iy} D_h(\vp). \]
This limit can be explicitly computed by an application of the strong Szeg\"o limit theorem.
(See, e.g.,~\cite{BS})
and thereby the Cauchy identity reappears.

The Toeplitz determinant,  or the Fredholm determinant
coming from the Borodin-Okounkov identity~\cite{BO, BW}, is
 the starting point in the analysis of limit laws
for $\la_1$.  This analysis together with the Robinson-Schensted-Knuth (RSK) correspondence
gives a new class of limit laws,
first discovered in the context of random matrix
theory~\cite{TW1, TW2},  for a number of probability models.
Indeed, the result of Baik, Deift and Johansson~\cite{BDJ} for the limit law
of the length, $\ell_N(\pi)$, of the longest increasing subsequence in a random permutation
$\pi\in \cS_N$
is the now classic example:  Exponential specialization  of the
Gessel identity
together with the RSK correspondence\footnote{RSK
 associates bijectively to each permutation $\pi$ a pair
of standard Young tableaux $(P,Q)$ of
the same shape $\la$ such that $\ell_N(\pi)=\la_1$. See, e.g.,~\cite{St}.}
shows that
\[ \sum_{N=0}^\iy \pr\left(\ell_N\le h\right)\,{t^N\ov N!} \]
is an $h\times h$ Toeplitz determinant with symbol
$\vp(z) =e^{\sqrt{t}(z+1/z)}$. An  asymptotic analysis
of this Toeplitz determinant (using the
steepest descent method for Riemann-Hilbert problems~\cite{DZ})
shows that
\[ \lim_{n\ra\iy}\pr\left({\ell_N -2\sqrt{N}\ov N^{1/6}}< s\right)=F_2(s) \]
where $F_2$ is the limiting distribution of the largest eigenvalue (suitably centered and normalized)
in the Gaussian Unitary Ensemble~\cite{TW1}.
Similar results hold for longest
increasing subsequences in symmetrized random permutations~\cite{BR1, BR2} and
random words~\cite{AvM, Jo2, ITW, TW3}, for height fluctuations
in various growth models~\cite{BR3, GTW1, GTW2, Jo1},
and tiling problems~\cite{Jo3}, as well as extensions
to the other rows of the Young diagram~\cite{BOO, Jo2,Ok2}.

In the theory of symmetric functions there are many important
generalizations of Schur functions~\cite{Mac}.  These generalizations
satisfy Cauchy identities and it is  natural
to inquire about more general Gessel identities. However, one quickly sees
that without determinantal formulas of the type that exist for Schur functions
(the Jacobi-Trudi identity), Gessel identities seem unlikely.  Nevertheless,
the question of possible limit laws for sums of the type (\ref{GesselLHS})
remains interesting.

This paper initiates work in this direction.  Instead of Schur functions we shall work with
Schur $Q$-functions which have pfaffian representations.
 These functions, introduced by Schur in 1911 in his analysis of the projective
representation of the symmetric group, now have a combinatorial theory that parallels the combinatorial
theory of Schur functions.
This theory, due to Sagan~\cite{Sa} and Worley~\cite{Wo}
(see also \cite{Ste, HH}), is based on a shifted
version of the RSK algorithm.   Whereas the usual RSK algorithm associates bijectively to each
$\bN$-matrix $A$ a pair of semistandard Young tableaux, the shifted RSK algorithm associates
bijectively to each $\bP$-matrix\footnote{Informally, a $\bP$-matrix is an $\bN$-matrix where we allow
the nonzero entries to be either marked or unmarked.  Precise definitions are given below.} $A$
 a pair of shifted Young tableaux.  There is a notion of increasing paths
and the length of the maximal path, $L(A)$,  equals the number of boxes in the first row of the shifted tableau.
Thus  it is natural to assign to each
partition $\la$ into distinct parts, i.e., a \textit{strict partition}, the probability
\be \pr\left(\{\la\}\right)= {1\ov Z}\, Q_\la(x) \, P_\la(y) \label{shiftedSchur} \ee
where $Q_\la$ and $P_\la$ are
the Schur  $Q$-functions and
$Z$ is a normalization constant.   We call this measure the \textit{shifted Schur measure}.

At first our analysis is for general parameters
$x$ and $y$ appearing in the shifted Schur measure, and we find that there is indeed
a Gessel identity.  (It follows from the Ishikawa-Wakayama pfaffian summation 
formula~\cite{IW}.)
Then we specialize the measure by
choosing the first $m$ coordinates of $x$ and the first $n$ coordinates of $y$ equal to
$\al$ ($0<\al<1$) and the rest equal to zero. We call this $\al$-specialization and denote
the resulting specialized shifted Schur measure
by $\pr_\sigma$ where $\sigma=(m,n,\alpha)$ denotes the parameters of the measure.
Now, however, the matrix on the right side of the Gessel identity is not
Toeplitz and so the earlier analytical methods are not immediately
available to us.  Nevertheless, we do find that the distribution function
for $L(A)=\la_1$ can be expressed in terms of the Fredholm determinant
of an operator which is a perturbation of a direct sum of products of Hankel operators.
In the end we can show that the trace norm of the perturbations tend to zero and are able
to determine the asymptotics.

We asume\footnote{The stated restriction on $\tau$ is
very likely unnecessary for the validity
of the final result. Some details of the proof would be different in the other cases
but we did not carry them out.} that  $m/n=\tau$  is a constant satisfying
$\al^2<\tau<\al^{-2}$.   Our main result is\sp

\noindent{\bf Theorem}
{\it Let $\pr_{\sigma}$ denote the  $\al$-specialized
shifted Schur measure with $\tau$ satisfying the stated restriction.
Then there exist constants $c_1=c_1(\al,\tau)$ and $c_2=c_2(\al,\tau)$ such that}
\[ \lim_{n\ra\infty}\pr_{\sigma}\left({\la_1-c_1\,n\ov c_2\, n^{1/3}}<s\right) = F_2(s). \]
\sp

{}For $\tau=1$ the constants have a particularly simple form; namely,
\[ c_1(\al,1)= {4\al\ov 1-\al^2}\>\>\> \textrm{and} \>\>\> c_2(\al,1)= {\left(2\al(1+6\al^2+\al^4)\right)^{1/3}\ov 1-\al^2}\> . \]
Expressions for $c_1$ and $c_2$ in general are given  in \S6.
For the Schur measure the corresponding $\al$-specialization limit law was derived by
Johansson \cite{Jo1}; as the
theorem shows, the two limit laws are identical.\footnote{We note that our $\al$ is related to Johansson's
$q$ by $q=\al^2$.  For the $\al$-specialized Schur measure, $c_1(\al,1)=2\al/(1-\al)$
and $c_2(\al,1)=\al^{1/3}(1+\al)^{1/3}/(1-\al)$.}
The table of contents provides a description of the layout of this paper.

\section{Schur $Q$-Functions}
\setcounter{equation}{0}
This section and the next summarize the properties of the Schur $Q$-functions
and the shifted RSK algorithm
that we will need in this paper. The material  is not new on our part. It
is presented
to establish the notation used in subsequent sections as well as a
convenience to the reader.
A complete presentations can be found in  the books by Macdonald~\cite{Mac},
Hoffman and Humphreys~\cite{HH}, the papers by
Sagan~\cite{Sa}, Stembridge~\cite{Ste}, and the
thesis of   Worley~\cite{Wo}.

If $\la=(\la_1,\la_2,\dots)$ is a partition of $n$, we denote this by $\la\vdash n$.
The length of $\la$  is denoted by $\ell(\la)$.
Let $\cP_n$ denote the set of all partitions of $n$ and $\cP:=\bigcup_{n=0}^\iy\cP_n$.
($\cP_0$ is the empty partition.) Introduce
$\cD_n\subset \cP_n$:  the set of partitions of $n$ into \textit{distinct parts}.
 For example
\[ \cD_6=\left\{\left\{6 \right\},\left\{5,1\right\},\left\{4,2\right\},
        \left\{3,2,1\right\}\right\}.\]
Let
$\cD:=\bigcup_{n=0}^\iy \cD_n $,
the set of all partitions into distinct parts.
If $\la\vdash n$ is a partition with distinct parts, we denote this by
$\la\models n$ and call $\la$ a \textit{strict} partition of $n$.

Associated to a strict partition $\la$ is a \textit{shifted shape} $S(\la)$.  One
starts with the usual Young diagram $Y(\la)$ and for $i=1,2,\ldots,\ell(\la)$ simply
indents the $i^{\scriptstyle\textrm{th}}$ row to the right
by $i-1$ boxes. The result is $S(\la)$.
We usually use  $\la\in\cD$ both
to denote a strict partition and the shifted shape $S(\la)$.

We let $\bN$ denote the set of positve integers and
\[ \bP=\left\{\1p, 1, \2p, 2, \3p, 3, \ldots\right\}\]
with the ordering
\[ \1p<1<\2p<2<\3p<3<\cdots. \]
We call the elements either \textit{marked} or \textit{unmarked} depending
on whether the element is primed or not.  When we do not wish to distinguish
a marked element $m^\prime$ from the unmarked element $m$, we write $m^*$.
A \textit{(generalized) shifted Young tableau}, $T$, is an assignment of
elements of $\bP$ to a shifted shape $\la$ having the properties
\begin{description}
\item[T1] $T$ is weakly increasing across rows and down columns.
\item[T2] For each integer $m^*$, there is at most one $m^\prime$ in
each row and at most one $m$ in each column of $T$.  (Thus the marked elements
are strictly increasing across rows
of $T$ and the unmarked elements are strictly
increasing down columns of $T$.)
\end{description}
An example of a shifted tableau of shape $(7,5,3,2,1)$ is
\[\begin{array}{ccccccc}
\1p&1&\2p&2&2&\5p&6\\
&\2p&2&\3p&4&5\\
&&\3p&4&5\\
&&&6&\7p\\
&&&&\7p.
\end{array}\]
To each shifted tableau $T$ we associate a monomial
\[ x^T=x_1^{a_1}x_2^{a_2}\cdots x_{m}^{a_m}\cdots\]
where $a_m$ is the number of times $m^*$ appears in $T$.  Thus as far
as the monomial is concerned, we do not distinguish between marked and unmarked
elements.  In the above example,
\[ x^T=x_1^2 x_2^5 x_3^2 x_4^2 x_5^3 x_6^2 x_7^2. \]

Let $\la$ be a strict partition of $n$.  The Schur $Q$-function,
the generating function of shifted tableaux of shape $\la$,  is
\be Q_\la(x):=\sum_{T} x^T, \label{Qfn}\ee
where the sum runs over all shifted tableaux of shape $\la\models n$.
  The Schur $Q$-function is the analogue of the Schur function
$s_\la$ when one replaces semistandard Young tableaux of shape $\la$ by
shifted tableaux of shape $\la$.  (Of course, here $\la$ must be a strict
partition.)
  It will be convenient to introduce the Schur $P$-function
\[ P_\la(x)={1\ov 2^{\ell(\la)}} \, Q_\la(x).\]

We remark that a shifted tableau $T$
of shifted shape $\la\models n$ is called \textit{standard} if it has
no marked elements and uses each unmarked letter $1,2,\ldots, n$ exactly once.
Schur showed that the number of standard shifted tableaux of shape
$\la\models n$, $\la=(\la_1,\la_2,\ldots,\la_\ell)$  is
\begin{equation}
 f_s^\la ={n!\ov \la_1! \la_2!\cdots \la_{\ell}!}\prod_{1\le i<j\le\ell}
{\la_i-\la_j\ov \la_i+\la_j}\, . \label{shiftedSYT}
\end{equation}
This should be compared with the number of standard Young tableaux
\[ f^\la={n!\ov \ell_1!\ell_2!\cdots\ell_k!}\prod_{1\le i<j\le k}(\ell_i-\ell_j)\]
where $\la=(\la_1,\la_2,\ldots,\la_k)$ and $\ell_j=\la_j+k-j$.
The number of semistandard
Young tableaux of shape $\la$ that can be formed using the integers $1,2,\ldots,n$
is
\[ d^\la(n)=s_\la(\overbrace{1,\ldots,1}^n,0,0,\ldots). \]
Similarly, the number of shifted tableaux of shape $\la$ that can be formed using
the integers $1^*, 2^*,\ldots, n^*$ is
\[ d^\la_s(n)=Q_\la(\overbrace{1,\ldots,1}^n,0,0,\ldots).\]
This specialization of $Q_\la$ will be important below.

The Schur $Q$-functions satisfy a Cauchy identity
\be \sum_{\la\in\cD} Q_\la(x) \, P_{\la}(y) =
\prod_{i, j=1}^\iy { 1+x_i y_j\ov 1-x_i y_j}=Z. \label{Qcauchy}\ee
The right-hand side enumerates all matrices $A$ whose entries are chosen
from $\bP\cup \{0\}$: the denominator counts matrices with entries
in $\bN\cup \{0\}$, while the numerator accounts for the primes.
We call these matrices  $\bP$-matrices.    The above
product will  frequently be specialized to $x=(x_1,x_2,\ldots,x_m,0,\ldots)$
and $y=(y_1,y_2,\ldots,y_n,0,\ldots)$.  We use the same symbol $Z$ to
denote this specialization.  It will be clear from the context how to
interpret $Z$.

Define symmetric functions $q_k$ by
\be
Q(t):=\prod_{i=1}^\iy{1+t x_i\ov 1-t x_i} = \sum_{k=0}^\iy q_k(x) t^k\, .
 \label{schurQseries}\ee
(When necessary to
indicate the dependence upon $x$,  we write $Q(t,x)$.)
It  follows from $Q(t)Q(-t)=1$   that
\[ q_{2m}=\sum_{r=1}^{m-1} (-1)^{r-1} q_r q_{2m-r}+\hf (-1)^m q_m^2,\]
which shows that $q_{2m}\in \bQ[q_1,q_2,\ldots,q_{2m-1}]$ and hence by induction
on $m$,
\[ q_{2m}\in\bQ[q_1,q_3,q_5,\ldots,q_{2m-1}]. \]
Denote by $\Gamma$ the subring of $\Lambda$ generated by the $q_r$:
\[ \Gamma=\bZ[q_1,q_3,\ldots].\]
If $\la=(\la_1,\la_2,\ldots)$ we let
\[ q_\la:=q_{\la_1}q_{\la_2}\cdots. \]
It is known that the $q_\la$ with $\la$ strict form a $\bZ$-basis of $\Gamma$.

We now give the classical definition of the Schur $Q$-function.
(Of course, in this presentation  it is a theorem.)
If $\la$ is a strict partition of length $\le n$,  then $Q_\la$ equals
the coefficient of $t^\la:=t_1^{\la_1}t_2^{\la_2}\cdots$ in
\[ Q(t_1,t_2,\ldots,t_n)=\prod_{i=1}^n Q(t_i)\prod_{i<j} F(t_i^{-1}t_j) \]
where
\[ F(y)={1-y\ov 1+y} = 1 + 2\sum_{r\ge 1}(-1)^r y^r \]
and $Q$ is defined by (\ref{schurQseries}); in particular, for $r>s$
\be Q_{(r,s)}=\left({1\ov 2\pi \mi}\right)^2\int\!\!\!\int t_1^{-r-1} t_2^{-s-1}
        F(t_2/t_1) Q(t_1) Q(t_2)\, \df t_1\, \df t_2 \label{Qrs}\ee
where the contours could be chosen to be circles with
$\vert t_2\vert<\vert t_1\vert$.\footnote{This requires that $x\in\ell^1$
and that the poles $x_i\inv$ lie outside the contours. Notice that if the $t_1$
contour were deformed to one inside the $t_2$ contour then since $Q(t)Q(-t)=1$ the residue
at the pole $t_1=-t_2$
crossed would be $2\,t_2^{-r-s-1}$. The integral of this equals zero as long as $r$
and $s$ are not both zero. This shows that the contours can also be chosen so that
$\vert t_2\vert>\vert t_1\vert$. Equivalently, the integral representation holds for
$r<s$ as well.}\sp

Here are some additional properties of Schur $Q$-functions:
\begin{enumerate}
\item The $Q_\la$, $\la$ strict, form a $\bZ$-basis of $\Gamma$.
\item Using (\ref{Qrs}) and the expansion for $F$ we have for $r>s\ge 0$
\be Q_{(r,s)}=q_r q_s + 2\sum_{i=1}^s (-1)^i q_{r+i}q_{s-i} \label{QqDefn}\ee
For $r\le s$ we define $Q_{(r,s)}=-Q_{(s,r)}$.  Now let $\la$ be a
strict partition which we write in the form $\la=(\la_1,\la_2,\ldots,\la_{2n})$
where $\la_1>\la_2>\cdots>\la_{2n}\ge 0$.  Define the $2n\times 2n$ antisymmetric
matrix
\[ M_\la=\left(Q_{(\la_i,\la_j})\right), \]
then we have
\be Q_\la= \pf(M_\la) \label{Qpfaffian} \ee
where \pf\ denotes the pfaffian.
\end{enumerate}

\section{Shifted RSK Algorithm\label{RSKsec}}
\setcounter{equation}{0}
{}For later convenience, we use a nonstandard labeling
of the matrix $A$: rows are numbered starting at the lower left-hand
corner of $A$ and columns have the usual left-to-right labeling.
To each  $\bP$-matrix $A$ we (bijectively) associate a biword $w_A$
as follows.
For a fixed column index we scan the matrix for increasing values
of the row index.  If the $(i,j)$ entry is unmarked with value
$a_{ij}$ we repeat the pair
$ \left(\begin{array}{c} j \\ i  \end{array}\right)$
$a_{ij}$ times in $w_A$.  If the
$(i,j)$ element is marked,  the $i$ of the first pair
$ \left(\begin{array}{c} j\\ i \end{array}\right) $
appearing in $w_A$ is marked.
For example, if
\[ A=\left(\begin{array}{ccc}
        1&2&0\\
        \1p&0&\2p\\
        \3p&0&1\end{array}\right)\]
then
\[ w_A=\left(\begin{array}{cccccccccc}
                1&  1&  1&  1&  1&  2&  2&  3&  3&  3\\
              \1p&  1&  1&\2p&  3&  3&  3&  1&\2p&  2
        \end{array}\right).\]

A description of the shifted RSK algorithm is
more involved than the usual RSK algorithm, though
the general features remain the same. Namely, there
is a row bumping (and column bumping) procedure which
when iterated on a sequence $\al$ whose elements are
in $\bP$ gives a shifted tableaux $S$, the insertion
tableaux.  (This is applied to the sequence in the bottom
half of the biword $w_A$.)  The top half of $w_A$ gives
a recording tableaux $T$.   We now state the final result, referring
the interested reader to either \cite{HH} or \cite{Sa}.\sp

\noindent{\bf Theorem} (Sagan, Worley).
{\it There is a bijective correspondence between $\bP$-matrices $A=(a_{ij})$
and ordered pairs $(S,T)$ of shifted tableaux of the same shape, such
that $T$ has no marked letters on its main diagonal.  The correspondence
has the property that $\sum_i a_{ij}^*$ is the number of entries $t$ of $T$
for which $t^*=j$, whereas $\sum_j a_{ij}^*$ is the number of entries $s$
of $S$ for which $s^*=i$.  (Recall that the $^*$ means we do not distinguish
between a marked or unmarked form of an integer.)}\sp

We call the matrix $S$, respectively $T$, to be of type $s$, respectively $t$.

An important property of the RSK algorithm is its relationship to
increasing subsequences of maximal length in the biword $w_A$ (equivalently,
increasing paths in the matrix $A$ of maximal weight).  The shifted RSK algorithm
of Sagan and Worley shares a similar property once the notion of an increasing
subsequence is properly formulated.
Let $\psi\in\bP^m$, $\phi\in\bP^n$, and denote by
$\psi\sqcup\phi$  the concatenation
$(\psi_1,\ldots,\psi_m,\phi_1,\ldots,\phi_n)$.  Denote by $\textrm{rev}(\phi)$
the reverse of $\phi$, $\textrm{rev}(\phi)=(\phi_n,\phi_{n-1},\ldots,\phi_1)$.
Given a sequence $\al$ from $\bP$, an \textit{ascent pair} $(\psi,\phi)$ for $\al$
is a pair of subsequences $\psi$ of $\textrm{rev}(\al)$ and $\phi$ of $\al$
such that if $\psi\in\bP^m$ and $\phi\in\bP^n$ then
\begin{enumerate}
\item $\psi\sqcup\phi$ is weakly increasing with respect to the ordering of $\bP$.
\item For all $k\in\bN$, at most one (unmarked) $k$ appears in $\psi$.
\item For all $k\in\bN$, at most one (marked) $k^\prime$ appears in $\phi$.
\end{enumerate}
Thus the unmarked symbols are strictly increasing in $\psi$ and the marked symbols
are strictly increasing in $\phi$.
The length of $\psi\sqcup\phi$ is defined to be $m+n-1$.
(Note that the length here is one less than the length defined in
either \cite{HH} or \cite{Sa}.)
Let $\al\in\bP^n$ and let $L(\al)$ denote the length of the longest ascent
pair $(\psi,\phi)$ of $\al$.  Then we have\sp

\noindent{\bf Theorem} (Sagan, Worley).
{\it If $\al$ is a sequence from $\bP$ and $T$ is the shifted tableaux
of shape $(\la_1,\la_2,\ldots,\la_k)$
resulting from the insertion of $\al$ (following the rules of
the shifted RSK algorithm), then
$L(\al)=\la_1$.}\sp

Here is an example of a increasing path  displayed in the $\bP$-matrix $A$.
In this example $L(A)=16$.\footnote{Note, for example, that the path segment
$\1p\ra\2p$ in the upper right hand corner contributes a weight of two, not three,
since the marked elements are strictly increasing.}
\[
\left(
\begin{array}{*{13}c}
\1p & & 0 & & \3p & &  0 & & \2p & &  \mathbf{\,\1p} & \ra &\mathbf{\2p} \\
0 & & \1p & &0 & &\mathbf{1}&\ra & 0 &\ra & \stackrel{\uparrow}{\mathbf{3}} & & 0\\
1 & & 2 & &  0 & & \stackrel{\uparrow}{\mathbf{1}} &  &0 && 0 & &\3p \\
1 &&0 && \1p && \stackrel{\uparrow}{\mathbf{\,\2p}} && \1p && 0 && 0 \\
0 && 1 && 0 && \stackrel{\uparrow}{\mathbf{\,\1p}} && 0 && 0 & &\1p \\
&&&&&\nearrow &&&&&&&\\
0 && \1p && \mathbf{\1p} && 0 && 0 && 0 && 0 \\
0 && 0 && \stackrel{\uparrow}{\mathbf{\1p}}&  \leftarrow&\mathbf{\1p}
 & \leftarrow&\mathbf{\1p} & \leftarrow&\mathbf{\1p}&\leftarrow&\mathbf{\2p}
\end{array}
\right)
\]
A quick way to compute this length is to apply a modified patience sorting
algorithm~\cite{AD} to the lower
row of the associated biword. (We leave this as an exercise to the interested reader.)\ \
Of course, another way is to apply the shifted RSK algorithm and count
the number of boxes in the first row.

\section{A Gessel Identity}
The Gessel identity~\cite{Ge} (see also \cite{TW3})  states that the sum
\[ \sum_{{\la\in\cP \atop \la_1\le h}} s_\la(x)\, s_\la(y) \]
equals a certain $h\times h$ Toeplitz determinant.  The proof begins
by expressing the Schur functions $s_\la$ as determinants (the
Jacobi-Trudi identity) and proceeds by recognizing this sum of
products of determinants
as the expansion of a single determinant
of the product of two (nonsquare) matrices. (This expansion is called
the Cauchy-Binet expansion.)

We are interested in sums of the form
\be \sum_{{\la\in\cD \atop \la_1\le h}} Q_\la(x)\, P_\la(y) \label{ourSum}\ee
where now, as we have seen, the $Q_\la$ and $P_\la$ are given by pfaffians.
What is needed is a pfaffian version of the Cauchy-Binet
formula.  Fortunately, Ishikawa
and Wakayama~\cite{IW} have such a  formula.  (See also, Stembridge~\cite{Ste}.)

Introduce
\[ I_r^h=\left\{I=(i_1,\ldots,i_r): 1\le i_1<\cdots<i_r\le h\right\}\]
and denote by $A_J$ the submatrix formed from $A$ by taking those rows and columns
indexed by $J\in I_r^h$. Then the pfaffian summation formula is\sp

\noindent{\bf Theorem}\ (Ishikawa-Wakayama).
{\it Let $A=\left(a_{ij}\right)_{0\le i,j \le h}$ and
$B=\left(b_{ij}\right)_{0\le i,j \le h}$ be $(h+1)\times (h+1)$ skew symmetric
matrices, $h\in\bN$. Then
\ba
\sum_{{0\le r\le h \atop r\>{ \textrm{\tiny even}}}} \sum_{I\in I_r^h} \ga^{|I|}
\pf\left(A_I\right)\,\pf\left(B_I\right) +
\sum_{{0\le r\le h \atop r\>{ \textrm{\tiny odd}}}} \sum_{I\in I_r^h} \ga^{|I|}
\pf\left(A_{0I}\right)\, \pf\left(B_{0I}\right)& \\
& \hspace{-15ex}=(-1)^{h(h-1)/2}\;
\pf\left(
\begin{array}{ll}
-A & I_{h+1}\\
-I_{h+1} & C
\end{array}\right)
\label{summation4}
\ea
where $C=(C_{ij})_{0\le i,j\le h}$ is the $(h+1)\times (h+1)$ antisymmetric matrix
\[ C_{ij}=\left\{\begin{array}{l}
                \ga^j \, b_{0j},\;\textrm{if}\;\; i=0, j\ge 1, \\
                \ga^i \, b_{i0},\;\textrm{if}\; i\ge 1, j=0,\\
                \ga^{i+j}\, b_{ij}, \;\textrm{if}\; i\ge 1, j\ge 1.
                \end{array}\right.\]
Here $|I|=\sum i_k$ and $A_{0I}$,  where $I=\{i_1,i_2,\ldots\}$, stands for $A_J$ where
$J=\{0,i_1,i_2,\ldots\}$.}\sp

Observe that $I_r^h$ is the set of partitions with exactly $r$ distinct parts such
that the largest part is less than or equal to $h$.  For any such partition,
$r\le h$. In  (\ref{ourSum}) we break
the sum into two sums---the first sum over distinct partitions with an even number of
parts and the second sum over distinct partitions with an odd number of parts.
Recalling the pfaffian representation (\ref{Qpfaffian}) of $Q_\la$, we note that
if $\la$ has an odd number of parts then we extend the partition by appending $0$
giving us a vector of even length.  Thus the sum appearing in
the pfaffian summation formula is (up to a reversal of labels)  the sum over
distinct partitions satisfying $\la_1\le h$.  From (\ref{Qpfaffian}) we see
that $A_{h}(x)$ is the $(h+1)\times(h+1)$ antisymmetric matrix
\[\left(\begin{array}{cc}
                0&-q^t\\
        q & \widehat{Q}_h(x) \end{array}\right)\, .\]
Here $q$ is the $h\times 1$ matrix with elements
$q_r(x)$ ($r=1,2,\ldots,h$), $\widehat{Q}_h(x)$ is the \linebreak $h\times h$ antisymmetric matrix
with elements $Q_{(r,s)}(x)$, and $q^t$ denotes the transpose of $q$.
(Recall that $Q_{(r,0)}=q_r$.)
The pfaffian representation for $P_\la$ is obtained from the $Q_\la$ pfaffian
representation by inserting the factor $2^{-\ell(\la)}$. The matrix $B_h$ is
\ba B_{h}(y)&=&\left(\begin{array}{cc}
                0& -{1\ov 2}\,q^t\,\vspace{1ex}\\
        {1\ov 2}\, q & {1\ov 4} \,\, \widehat{Q}_h(y) \,
\end{array}\right)\vspace{2ex}\\
&=&\left(\begin{array}{cc}
                1 & 0\\
               0 & {1\ov 2}
                \end{array}\right)\;
                A_{h}(y)\,
                \left(\begin{array}{cc}
                1 & 0\\
               0 & {1\ov 2}\end{array} \right).\ea
Applying the summation formula
then gives for $h\in\bN$
\ba \sum_{{\la\in\cD \atop \la_1\le h}} \,Q_\la(x) \, P_\la(y) &=&
\pf\left(
\begin{array}{ll}
-A_{h}(x)& I \\
-I & B_{h}(y)
\end{array}\right)\\
&=& \left(\det\left(I-A_{h}(x) B_{h}(y)\right)\right)^{1/2}.
\ea
The $\pm$ factors are accounted for by the reversal of labels; or
more simply,  because we are computing probabilities.

Introducing the antisymmetric matrix
\[ K_{h}(x)=\left(\begin{array}{cc}
                1 & 0\\
               0 & {1\ov \sqrt{2}}
                \end{array}\right)\:
        A_{h}(x)\:\left(\begin{array}{cc}
                1 & 0\\
               0 & {1\ov \sqrt{2}}
                \end{array}\right), \]
or more explicitly,
\be K_{h}(x)_{rs}=\left\{\begin{array}{lr}
                        -{1\ov\sqrt{2}}\,q_s(x) &\:\: r=0,\: s\ge 1,\vspace{1ex} \\
                        {1\ov\sqrt{2}}\,q_r(x) &\:\:  r\ge 1, \: s=0,\vspace{1ex} \\
                        {1\ov 2}\, \, Q_{(r,s)}(x) \, &\:\: r\ge 1, s\ge 1,
                        \end{array}\right.\label{Kmatrix}\ee
we obtain our Gessel identity:
\be
\sum_{{\la\in\cD \atop \la_1\le h}} \,Q_\la(x) \, P_\la(y) =
\Big(\det\left(I-K_{h}(x)K_{h}(y)\right)\Big)^{1/2}\, .\label{Gessel}
\ee
Observe that it
follows from this and the Cauchy identity (\ref{Qcauchy}) that
\be \lim_{h\ra\iy}\det\left(I-K_{h}(x) K_{h}(y)\right)=Z^2.\label{Qcauchy2}\ee

\section{Shifted Schur Measure}
\setcounter{equation}{0}
Let $\cP_{m,n}$ denote
the set of $\bP$-matrices of size $m\times n$. For $A\in\cP_{m,n}$
we recall that  $L(A)$  denotes the
length of the longest increasing path in  $A$.
Let $x=(x_1,x_2,\ldots)$ and $y=(y_1,y_2,\ldots)$
with $0\le x_i<1$ and $0\le y_i<1$.
We assume the matrix elements $a_{ij}$ are
distributed independently with a geometric distribution with parameter
$x_i\,y_j$.
Specifically, for $k\ge 1$
\[\pr\left(a_{ij}=k\right)=
\pr\left(a_{ij}=k^\prime\right)
= \left({1-x_i y_j\ov 1+ x_i y_j}\right) \,(x_i y_j)^k \]
and
\[ \pr\left(a_{ij}=0\right)={1-x_i y_j\ov 1+x_i y_j}\, .\]
We have, of course,
\[ \sum_{{k^{*}}=0}^\iy \pr\left(a_{ij}=k^*\right)={1-x_i y_j \ov 1+x_i y_j}
+2\sum_{k\ge 1} \left({1-x_i y_j\ov 1+ x_i y_j}\right)(x_i y_j)^k=1\]

Let $\mathcal{P}_{m,n,s,t}\,$
($t\in \bN^m$, $s\in \bN^n$) denote the set of
 $A\in\bP_{m,n}$ satisfying, for $1\le i\le m$ and $1\le j\le n$,
\[  \sum_{1\le j\le n} a_{ij}^*=s_i
\;\textrm{and}\;\sum_{1\le i\le m} a_{ij}^* = t_j\]
Then for $A\in\mathcal{P}_{m,n,s,t}$ we have
\[ \pr\left(\{A\}\right)= \prod_{{1\le i \le m \atop 1\le j \le n}}
\left(1-x_i y_j\ov 1+ x_i y_j\right)\, x^s y^t\,= {1\ov Z} \, x^s y^t . \]
Since right hand side does not depend upon the
$A$ chosen in $\mathcal{P}_{m,n,s,t}$,   the conditional probability
\[ \pr\left(L\le h \mid \sum_j a_{ij}^*=s_i, \sum_i a_{ij}^*=t_j\right)\]
is uniform.
Note that this  uses both the independence and the geometric
distribution of the random variables $a_{ij}$.

By the shifted RSK correspondence,
to each $A\in \mathcal{P}_{m,n,s,t}$ we associate bijectively
a pair $(S,T)$ of shifted tableaux of the same shape $\la\models N$
($N:=\sum_{i,j} a_{ij}^*$),
of types $s$ and $t$, respectively.
The condition $L(A)\le h$ becomes $\la_1\le h$.
Hence
\ba
\pr_{m,n}\left(L\le h\right)&=& \sum_{A\in \cP_{m,n}/\cP_{m,n,s,t}}
\pr\left(\ell(A)\le h | A\in \cP_{m,n,s,t}\right)\pr\left(A\in \cP_{m,n,s,t}\right)\\
        &=&\sum_{A\in \cP_{m,n}/\cP_{m,n,s,t}}\,
        {1\ov \vert \cP_{m,n,s,t}\vert} \,
{1\ov Z}\, x^s y^t \vert\cP_{m,n,s,t}\vert\\
&=&{1\ov Z}\,
\sum_{N\ge 0}\sum_{{\la\models N\atop \la_1\le h}} Q_\la(x) P_\la(y).\ea
(Here $\pr_{m,n}$ denotes probability before $\al$-specialization.) Thus, by (\ref{Gessel}),
\be \pr_{m,n}\left(L\le h\right)={1\ov Z}\,
\Big(\det\left(I-K_{h}(x) K_{h}(y)\right)\Big)^{1/2}.
\label{probrep}\ee
The above  used the combinatorial definition
(\ref{Qfn}) of the Schur $Q$-function.
The reason for the occurence of $P_\la(y)$ (instead of $Q_\la(y)$) is
that the recording tableaux $T$ has no marked elements on the diagonal
which accounts for the factor $2^{-\ell(\la)}$.  (There are exactly
$2^{\ell(\la)}$ entries on the main diagonal in a marked shifted tableaux
of shape $\la$ and so there are $2^{\ell(\la)}$ ways to mark and unmark
the diagonal elements.)

Observe the consequence that the distribution function $\pr_{m,n}\left(L\le h\right)$ is
a symmetric function of $x=(x_1,\ldots,x_m)$ and of $y=(y_1,\ldots,y_n)$.

\section{Proof of the Theorem}
\subsection{An Operator Formulation}
We begin  by deriving an alternative representation
for $\det\left(I-K_h(x) K_h(y)\right)$
in terms of Toeplitz and Hankel operators on the Hilbert space
$\ell^2(\bZ_{+})$ ($\bZ_+:=\bN\cup{0}$).
These may well be of independent interest in the theory of Schur $Q$-functions.
We assume at first only that $x,y\in\ell^1$ together with uniform estimates
$x_j,\,y_j\le c <1$

To set notation, we let $\left\{ e_j\right\}_{j\ge 0}$ denote
the canonical basis of $\ell^2(\bZ_{+})$.
Since the vector $e_0$ will occur frequently, we
denote $e_0$ by $e$, and it is convenient to set $e_{-1}=0$.
The backward shift operator $\La$ is characterized by
\[ \La e_j = e_{j-1} \]
and its adjoint $\La^*$ is the forward shift operator.  The two satisfy
\be \La\,\La^* = I \ \ \textrm{and}\ \
 \La^*\,\La = I - e\ot e, \label{LambdaIdentities}\ee
where for vectors $u$ and $v$ we denote by $u\ot v$ the operator sending a vector
$f$ to $u\,(v,f)$.

Suppressing temporarily the parameters $x$ and $y$ we define $L$ to be the matrix
with entries
\be L_{j\,k}=\left({1\ov 2\pi \mi}\right)^2
\intint {Q(z)\,Q(\z)\ov z^{j+1}\,\z^{k+1}}\, {\df z\,\df \z\ov z+\z}\> , \label{L}\ee
where $Q$ is defined in (\ref{schurQseries}).
Here the contours can be taken to be concentric circles of different radii near the unit
circle. Since $Q(z)\,Q(-z)=1$ the residue at $\z=-z$ in the integral defining
$L$ is zero, so we may freely choose whether the $z$-contour lies inside or outside of
the $\z$-contour without affecting the value of the integral. (Equivalently, $L$
is symmetric.)

{}From (\ref{Qrs})  we see that
\[ Q_{(j,k)}= L_{j-1,k} - L_{j,k-1}, \]
where we set $L_{-1,k}=0$. Note that $L_{j-1,0}=q_j$. The matrix elements $K_{j k}$ given
in (\ref{Kmatrix}) are then
\ba
K_{j\,k}=\left\{\begin{array}{ll}
        -{1\ov\sqrt{2}}\, L_{0,k-1},& j=0,k\ne 0,\vspace{1ex} \\
        \hspace{1.9ex}{1\ov \sqrt{2}}\, L_{j-1,0},& j\ne 0, k=0,\vspace{1ex}\\
        \hspace{2.5ex}{1\ov 2}\left(L_{j-1,k}-L_{j,k-1}\right), & j,k>0.\vspace{1ex}
        \end{array}\right.\ea
Introducing the vector $q=q(x)=(q_0(x),q_1(x),\ldots)$, the operator $K$ can 
be written
\be K={1\ov 2}\left(\La^* L - L\La\right) + {1\ov 2}\,\om\left(e\ot q - q\ot e\right),
\label{Koperator} \ee
where $\om=1-\sqrt{2}$.

The operator $L$ is expressible in terms of Toeplitz and Hankel matrices acting
on $\ell^2(\bZ_{+})$. Recall that $T(\psi)$, the Toeplitz matrix with symbol $\psi$,
has $j,k$ entry $\psi_{j-k}$ (subscripts denote Fourier coefficients here) while
the Hankel matrix $H(\psi)$ has $j,k$ entry $\psi_{j+k+1}$. If we assume
the contours in (\ref{L}) chosen so that $|\z|<|z|$ and expand $(z+\z)\inv$ in powers of
$\z/z$ we obtain
\[L_{i\,j}=\sum_{k=0}^\iy(-1)^k\left({1\ov 2\pi \mi}\right)^2\int\!\!\int z^{-i-k-2}\,
\z^{k-j-1}\,Q(z)\,Q(\z)\,\df z\,\df\z.\]
Now make the substitution $\z\to-\z\inv$ to obtain
\[(-1)^j\sum_{k=0}^\iy\left({1\ov 2\pi \mi}\right)^2\int\!\!\int z^{-i-k-2}\,\z^{-k+j-1}\,
{Q(z)\ov \tQ(\z)}\,\df z\,\df\z,\]
where $\tQ(\z)=Q(\z\inv)$. The $z$-integral gives $Q_{i+k+1}$
while the $\z$-integral gives $(\tQ\inv)_{k-j}$. It follows that
\be L=H(Q)\,T(\tQ\inv)\,J,\label{L1}\ee
where $J$ is the diagonal matrix with diagonal entries $(-1)^j$. If in the last
integrals we make the substitutions $z\to -z,\ \ \z\to-\z$ we find that also
\be L=-J\,H(Q\inv)\,T(\tQ).\label{L2}\ee

If we reintroduce our parameters $x$ and $y$, which we now write for notational
convenience as subscripts, and use the two representations of $L$
we see that
\[\det\,(I+L_x\,L_y)=\det\Big(I-H(Q_x)\,T(\tQ_x\inv)\,H(Q_y\inv)\,T(\tQ_y)\Big)\]
\[=\det\,\Big(I-T(\tQ_y)\,H(Q_x)\,T(\tQ_x\inv)\,H(Q_y\inv)\Big).\]
Here we have used the general identity $\det\,(I-AB)=\det\,(I-BA)$, valid if one
of the operators is trace class and the other bounded, and the fact that the Hankel
operators are Hilbert-Schmidt under our assumtions on $x$ and $y$. Another general fact is
\be T(\psi_1)\,H(\psi_2)+H(\psi_1)\,T(\widetilde\psi_2)=H(\psi_1\,\psi_2).
\label{THidentity}\ee
In particular, if $\psi_1$ is a minus function
(Fourier coefficients with positive index all vanish), then
$T(\psi_1)\,H(\psi_2)=H(\psi_1\psi_2)$.  {}From this we find that
\[T(\tQ_y)\,H(Q_x)=H(Q_x\,\tQ_y) \And T(\tQ_x\inv)\,H(Q_y\inv)=H(\tQ_x\inv\,Q_y\inv),\]
so the product of these equals $H(\phi)\,H(\widetilde\phi\inv)$, where
\be \phi(z):= Q_x(z)\,\tQ_y(z). \label{symbolPhi} \ee
Since yet another general identity is
\be T(\psi_1)\,T(\psi_2)=T(\psi_1\,\psi_2)-H(\psi_1)\,H(\widetilde\psi_2)\label{TTHH}\ee
we have $I-H(\phi)\,H(\widetilde\phi\inv)=T(\phi)\,T(\phi\inv)$, and we have shown that
\be\det\,(I+L_x\,L_y)=\det T(\phi) \,T(\phi^{-1}).\label{Ldet}\ee

If a symbol $\phi$ has geometric mean 1 and is sufficiently well behaved then
the strong Szeg\"o limit theorem says that
\[ \lim_{h\ra\iy} \det T_h(\phi) =E(\phi):=\exp\left(\sum_{n=1}^\iy n\, (\log\phi)_n \,
(\log\phi)_{-n}\right), \]
where $T_h(\phi)=(\phi_{j-k})_{j,k=0,\ldots,h-1}$. In the case of our symbol given by
(\ref{symbolPhi}) we find that
\[E(\phi)=\left(\prod_{i, j}{1+x_i\,y_j\ov 1-x_i\,y_j}\right)^2=Z^2,\]
where $Z$ is as in the right side of  (\ref{Qcauchy}). But there is another formula
for $E(\phi)$, namely~\cite{Wi}
\[ E(\phi)=\det T(\phi) T(\phi^{-1}),\]
and so from (\ref{Ldet}) we have the identity
\[\det(I+L_x\,L_y)=Z^2.\]

In the case of Schur functions, the right hand side of the Gessel identity
is a Toeplitz determinant, and the Cauchy identity for Schur functions emerges
as a consequence of the Szeg\"o limit theorem. In view of the last identity
and the connection between the operators $L$ and $K$
on $\ell^2(\bZ_{+})$ it is tempting to try to find, using these, an independent
derivation of (\ref{Qcauchy}). It will follow from (\ref{Koperator}) that
$I-K(x)\, K(y)$ and $I+L_x\,L_y$
differ by a finite rank operator. This operator cannot contribute to
the determinant but we do not see, {\it a priori}, why this is so. So such an
independent derivation eludes us.

To continue now, we let $P_h$ be the projection operator
onto the subspace of $\ell^2(\bZ_{+})$ spanned by $\{e_0,e_1,\ldots,e_h\}$.
Thus, if $K$ is the operator on $\ell^2(\bZ_{+})$ then
$K_h=P_h K P_h$.
Instead of working directly with the product $K(x)\,K(y)$, it will be convenient
to write $2\times 2$ matrices with operator entries.  Thus $\det\,(I-K_h(x)\,K_h(y))$
is the determinant of
\[\left(\begin{array}{cc} P_h&0\\0&P_h\end{array}\right)\,
\left(\begin{array}{cc} I &  K(x) \\
                  K(y) & I \end{array} \right)\,
\left(\begin{array}{cc} P_h&0\\0&P_h\end{array}\right)\, , \]
thought of as acting
on $P_h\,\ell^2(\bZ_+)\oplus P_h\,\ell^2(\bZ_+)$.
To simplify notation, we use $P_h$ to denote also
\[ \left(\begin{array}{cc} P_h&0\\0&P_h\end{array}\right)\,, \]
and set
\[ \cK=\left(\begin{array}{cc} 0 & K(x) \\ K(y) & 0\end{array}\right).\]
Thus,
\[\det\,(I-K_h(x)\,K_h(y))=\det\,P_h(I+\cK)P_h.\]

It follows from (\ref{Qcauchy2}) and the
infinite-dimensional version of Jacobi's theorem on the
principal $n \times n$ minor of the inverse of a (finite) matrix\footnote{The 
infinite-dimensional result follows by replacing the operator $\cK$ by
$P_N\cK P_N$, applying the finite-dimensional result and taking the $N\to\iy$ limit.
We use the fact that $\cK$ is trace class, which holds since the Hankel
operators are trace class. Of course all this requires that the infinite-dimensional operator
be invertible. This will follow from the limit results we establish, as we shall see
at the end of \S6.3.} that this may be written
\be\det\,(I-K_h(x)\,K_h(y))=Z^2\,\det\Big((I-P_h)(I+\cK)\inv (I-P_h)\Big).\label{Kdetrep}\ee
Thus our first goal is to compute $(I+\cK)\inv$.

Using the easily verified fact
\be \La^*\, L + L\,\La = q\ot q -e\ot e\label{LLambda}\ee
and (\ref{Koperator})
we see that $K=\La^* L +R^-=-L\La+R^+$
where
\be R^{\pm}={1\ov 2}\left(\pm q\ot q \mp e\ot e+ \om (e\ot q - q\ot e)\right).\label{R}\ee
Thus
\[ I+\cK=\twotwo{I}{\La^* L_x+R^-_x}{\La^* L_y+R^-_y}{I},\]
where subscripts have the usual meaning.

\subsection{Calculation of $\mathbf{(I+\cK)^{-1}}$}

The fundamental objects which will appear here are $I+L_x L_y$ and
$I+L_y L_x$ and we begin by showing that they are invertible and computing
their inverses. Define
\[H_1=H(Q_x\,\tQ_y\inv)\ \ \textrm{and}\ \  H_2=H(Q_y\,\tQ_x\inv).\]
We shall prove the basic identities
\be \left(I+L_x L_y\right)^{-1} = I - H_1 H_2\ \ \textrm{and}\ \
\left(I+L_y L_x\right)^{-1} = I-H_2 H_1.
\label{basicInverse}\ee
We do this with the help of (\ref{THidentity}) and (\ref{TTHH}). Using these and
(\ref{L1}) we find that
\[H_1\,H_2\,L_x=H(Q_x\,\tQ_y\inv)\,H(Q_y\,\tQ_x\inv)\,
H(Q_x)\,T(\tQ_x\inv)\,J\]
\[=H(Q_x\,\tQ_y\inv)\,\left[T(Q_y)-T(Q_y\,\tQ_x\inv)\,T(\tQ_x)\right]
\,T(\tQ_x\inv)\,J\]
\[=\left[H(Q_x)\,T(\tQ_x\inv)-
H(Q_x\,\tQ_y\inv)\,T(Q_y\,\tQ_x\inv)\right]\,J
=L_x-H_1\,T(Q_y\,\tQ_x\inv)\,J.\]
Thus
\be(I-H_1\,H_2)\,L_x=H_1\,T(Q_y\,\tQ_x\inv)\,J.\label{HHL}\ee
{}From this and (\ref{L2}) we obtain similarly
\[(I-H_1\,H_2)\,L_x\,L_y=-H_1\,T(Q_y\,\tQ_x\inv)\,H(Q_y\inv)\,T(\tQ_y)\]
\[=H_1\,H(Q_y\,\tQ_x\inv)\,T(\tQ_y\inv)\,T(\tQ_y)=H_1\,H_2.\]
This establishes the first identity of (\ref{basicInverse}), and the second is
obtained by interchanging $x$ and~$y$.

Beginning the calculation of $(I+\cK)\inv$ we refer to (\ref{LLambda}) and (\ref{R}) and
find that
\[\twotwo{I}{L_x \La}{L_y \La}{I}\,(I+\cK)=\twotwo{I+L_xL_y}{0}{0}{I+L_yL_x}
+\twotwo{L_x\La R^-_y}{R^+_x}{R^+_y}{L_y\La R^-_x},\]
where we use subscripts as before.
Using (\ref{basicInverse}) this may be written
\[\twotwo{I+L_xL_y}{0}{0}{I+L_yL_x}\left[I+\twotwo{(I-H_1H_2)L_x\La R^-_y}
{(I-H_1H_2)R^+_x}{(I-H_2H_1)R^+_y}{(I-H_2H_1)L_y\La R^-_x}\right].\]
Hence
\[(I+\cK)\inv=\left[I+\twotwo{(I-H_1H_2)L_x\La R^-_y}
{(I-H_1H_2)R^+_x}{(I-H_2H_1)R^+_y}{(I-H_2H_1)L_y\La R^-_x}\right]\inv\]
\[\times\twotwo{I-H_1H_2}{(I-H_1H_2)L_x\La}{(I-H_2H_1)L_y\La}{I-H_2H_1}.\]

To compute the entries of the matrix inside the large bracket we show that
\be(I-H_1H_2)L_x\La q_y=H_1\,H_2\,e\ \ \textrm{and}\ \ (I-H_1H_2)\,q_x=T_1\,e,
\label{Hqidentities}\ee
where we set
\[T_1=T(Q_x\,\tQ_y\inv)\ \ \textrm{and}\ \ T_2=T(Q_y\,\tQ_x\inv).\]
For the first we use the fact that $\La q=Le$, which gives
\[(I-H_1H_2)L_x\,\La q_y=(I-H_1\,H_2)\,L_x\,L_y\,e=H_1\,H_2\,e.\]
To derive the second we use the fact that $q_x=T(Q_x)\,e$ and compute
\[H_1\,H_2\,T(Q_x)=H(Q_x\,\tQ_y\inv)\,H(Q_y\,\tQ_x\inv)\,T(Q_x)
=H(Q_x\,\tQ_y\inv)\,H(Q_y)\]
\[=T(Q_x)-T(Q_x\,\tQ_y\inv)\,T(\tQ_y).\]
Since $T(\tQ_y)\,e=e$ (the matrix is upper-triangular with $0,0$ entry 1) this gives
\[H_1\,H_2\,q_x=q_x-T_1\,e,\]
which is equivalent to the desired identity. Of course the same identities hold
if we make the interchanges $x\leftrightarrow y$ and $1\leftrightarrow 2$.

With these identities and the fact that
$\La R^-=-{1\ov2}\La q\ot(q+\om e)$ we find that the matrix in large brackets
may be written
\[I+{1\ov2}\twotwo{-H_1H_2e\ot(q_y+\om e)}{\hspace{-8ex} T_1e\ot(q_x-\om e)+(I-H_1H_2)e\ot(\om q_x-e)}
{T_2e\ot(q_y-\om e)+(I-H_2H_1)e\ot (\om q_y-e)}{\hspace{-8ex}-H_2H_1e\ot(q_x+\om e)}.\]
This in turn has the form
\[I+\sum_{i=1}^4 a_i\ot b_i,\]
where
\[a_1={1\ov2}\twoone{-H_1H_2e}{-\om H_2H_1e+T_2e+\om e},\ \
a_2={1\ov2}\twoone{-\om H_1H_2e+T_1e+\om e}{-H_2H_1e},\]
\[a_3={1\ov2}\twoone{-\om H_1H_2e}{H_2H_1e-\om T_2e-e},\ \
a_4={1\ov2}\twoone{H_1H_2e-\om T_1e-e}{-\om H_2H_1e},\]
\[ b_1=\twoone{q_y}{0},\ \ b_2=\twoone{0}{q_x},\ \ \
b_3=\twoone{e}{0},\ \ b_4=\twoone{0}{e}.\]
At this stage we have shown that
\be (I+\cK)\inv=\left(I+\sum_{i=1}^4 a_i\ot b_i\right)\inv
\twotwo{I-H_1H_2}{(I-H_1H_2)L_x\La}{(I-H_2H_1)L_y\La}{I-H_2H_1}.\label{IplusKinv}\ee

If we have a finite rank operator
$\sum a_i\ot b_i$, then
\be \left(I+\sum a_i\ot b_i\right)^{-1}=I-\sum_{i,j} (S^{-1})_{ij}\, a_i\ot b_j,
\label{finiteRankInv}\ee
where $S$ is the matrix with entries
\[ S_{ij}=\delta_{ij}+(b_i,a_j).\]
In our case we have to compute 16 inner products, which is not as bad as it might seem since there
are basic inner products from which the others can be derived. And if we have evaluated
any inner product then we have evaluated another with the interchanges $x\leftrightarrow y$
and $1\leftrightarrow2$. Two basic inner products are trivial:
\[(e,e)=1 \ \ \textrm{and}\ \  (e,q_y)=1.\]
Two are not evaluable in simpler terms but just notationally. We set
\[t=(T_1e,e)=(T_2e,e) \And h=(H_1H_2e,e)=(H_2H_1e,e).\]
(The equality of the first two inner products follows from the facts that
$T_1^*=JT_2J$ and $Je=e$.) The nontrivial ones are
\[(T_1e,q_y)=1 \ \ \textrm{and}\ \  (H_1H_2e,q_y)=1-t.\]
For the first, we have $(T_1e,q_y)=(T_1e,T(Q_y)e)=(T(\tQ_y)T_1e,e)$,
and this is the $0,0$ entry of $T(\tQ_y)\,T(Q_x\tQ_y\inv)=T(Q_x)$. The $0,0$ entry equals 1.
For the second, we have
\[(H_1H_2e,q_y)=1-((I-H_1H_2)e,q_y)=1-(e,(I-H_2H_1)q_y)=1-(e,T_2e)\]
by the second identity of (\ref{Hqidentities}).

We can now write down all 16 inner products. For convenience we multiply them by $2$:
\[2\,(a_1,b_1)=-1+t,\ \ \ 2\,(a_2,b_1)=\om t+1,\ \ \ 2\,(a_3,b_1)=-\om (1-t),\ \ \
2\,(a_4,b_1)=-t-\om-2,\]
\[2\,(a_1,b_2)=\om t+1,\ \ \ 2\,(a_2,b_2)=-1+t,\ \ \ 2\,(a_3,b_2)=-t-\om-2,\ \ \
2\,(a_4,b_2)=-\om (1-t),\]
\[2\,(a_1,b_3)=-h,\ \ \ 2\,(a_2,b_3)=-\om h+t+\om, \ \ \ 2\,(a_3,b_3)=-\om h,\ \ \
2\,(a_4,b_3)=h-\om t-1,\]
\[2\,(a_1,b_4)=-\om h+t+\om,\ \ \ 2\,(a_2,b_4)=-h,\ \ \ 2\,(a_3,b_4)=h-\om t-1,
\ \ \ 2\,(a_4,b_4)=-\om h.\]

Let us see which vectors arise in the end. From (\ref{IplusKinv}) and
(\ref{finiteRankInv})
\be (I+\cK)\inv=\twotwo{I-H_1H_2}{(I-H_1H_2)L_x\La}{(I-H_2H_1)L_y\La}{I-H_2H_1}
-\sum_{i,j=1}^4 s_{ij}\,a_i\ot b_j',\label{Kinvrep}\ee
where $s_{ij}=(S\inv)_{ij}$ and
\[b_j'=\twotwo{I-H_2H_1}{\La^*L_y(I-H_1H_2)}{\La^*L_x(I-H_2H_1)}{I-H_1H_2}b_j\,.\]
The quantities that appear in the $a_i$, other
than $e$ which we can ignore since $(I-P_h)e=0$, are $H_1H_2e$ and $T_1e$ in the first
component and $H_2H_1e$ and $T_2e$ in the second.
Those in the $b_j$ are $q_y$ and $e$ in the first component and $q_x$ and $e$
in the second.
For $b_j'$ we use (\ref{Hqidentities}) to see that
$T_2e,\, H_2H_1e,\, \La^*L_yT_1e$ and $\La^*L_y(I-H_1H_2)e$ appear in the first
component and $T_1e,\, H_1H_2e,\, \La^*L_xT_2e$ and $\La^*L_x(I-H_2H_1)e$
in the second. (The $e$ which appear once again drop out in the end.)

Two new vectors appear here (as well as those obtained by the usual interchanges). We
claim that
\be \La^*L_yT_1e=T_2e-te\And \La^*L_y(I-H_1H_2)e=H_2H_1e+t\,T_2e-(1-h)e.\label{LaLy}\ee

{}For the first, we have
\[L_yT_1=-JH(Q_y\inv)T(\tQ_y)T(Q_x\tQ_y\inv)=-JH(Q_y\inv)T(Q_x)\]
\[=-JH(\tQ_xQ_y\inv)
=H(Q_y\tQ_x\inv)J,\]
so
\[\La^*L_yT_1e=\La^*H(Q_y\tQ_x\inv)e=T_2e-te.\]
For the second, we take transposes and interchange $x$ and $y$ in (\ref{HHL})
to obtain
\[L_y(1-H_1H_2)e=JT(\tQ_xQ_y\inv)H(Q_y\tQ_x\inv)e\]
\[=-T(Q_y\tQ_x\inv)H(\tQ_xQ_y\inv)e
=H(Q_y\tQ_x\inv)T(Q_x\tQ_y\inv)e.\]
So
\[\La^*L_y(1-H_1H_2)e=[H(zQ_y\tQ_x\inv)-e\ot T_2e]T(Q_x\tQ_y\inv)e\]
\[=H(zQ_y\tQ_x\inv)H(zQ_x\tQ_y\inv)e-(T_1e,T_2e)e
=H_2H_1e+(T_2e\ot T_1e)e-(T_1e,T_2e)e.\]
The next to last term equals $t\,T_2e$ while the last inner product equals
\[(JT_2^*JJT_1Je,e)=(T(Q_x\tQ_y\inv)T(\tQ_yQ_x\inv)e,e)=((I-H_1H_2)e,e).\]
This establishes the second claim.

It follows from the above that the only vectors that arise in the $b_j'$
are $T_2e$ and $H_2H_1e$ in the first component and $T_1e$ and $H_1H_2e$
in the second.

\subsection{Specialization}

\begin{figure}
\bc
\resizebox{7cm}{7cm}{\includegraphics{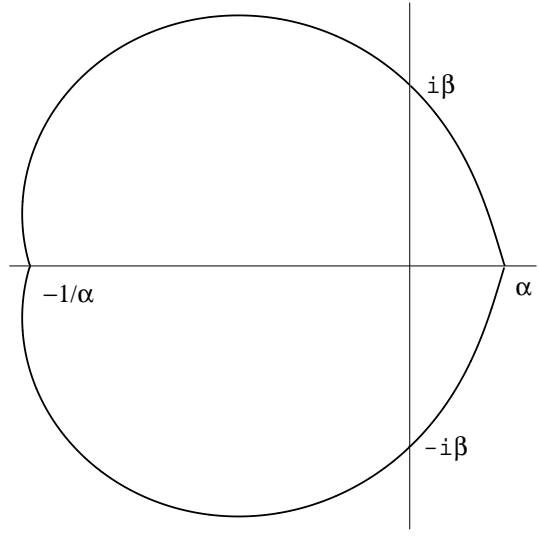}}
\caption{Steepest descent curve.}
\label{sd2Fig}
\ec
\end{figure}

At this point we impose the $\al$-specialization. Thus the first $m$ $x_i$
and the first $n$ $y_i$ are equal to $\al$ and the rest equal to zero. We assume
that $\tau=m/n$ is a constant satisfying
$\al^2<\tau<\al^{-2}$
and we first
determine the asymptotics as $n\to\iy$ of the quantities appearing in the inner
products. We claim
\[\lim_{n\ra\iy}t=0 \>\>\>\textrm{and}\>\>\> \lim_{n\ra\iy}h={1\ov2}.\]
For the first, we have
\be t=(T_1)_{0,0}={1\ov 2\pi \mi}\int\left({1+\al z\ov 1-\al z}\right)^m
\left({z-\al \ov z+\al }\right)^n{\df z\ov z}\,.\label{tintegral}\ee
If we apply steepest descent we see that the saddle points are to satisfy
\[{\tau\ov 1-\al^2z^2}+{1\ov z^2-\al^2}=0,\]
and so are given by
\[z=\pm \,\mi\,\sqrt{{1-\al^2\,\tau\ov\tau-\al^2}}=\pm \mi\beta.\]
These are purely imaginary under our assumption on $\tau$. The steepest descent curve
passes through these points and
closes at $\al$ and $-\al\inv$.   (See Fig.~\ref{sd2Fig}.)  The integral is $O(n^{-1/2})$.
{}For the second, we have
\[(H_1H_2)_{i,j}=\left({1\ov 2\pi \mi}\right)^2\sum_{k=0}^\iy\intint
\left({1+\al z\ov 1-\al z}\right)^m\left({z-\al \ov z+\al }\right)^n
\left({1+\al \z\ov 1-\al \z}\right)^n\left({\z-\al \ov \z+\al }\right)^m\]
\[\times z^{-i-k-2}\,\z^{-k-j-2}\,\df z\,\df\z\]
\vspace{.1ex}
\be=\left({1\ov 2\pi \mi}\right)^2\intint
\left({1+\al z\ov 1-\al z}\right)^m\left({z-\al \ov z+\al }\right)^n
\left({1+\al \z\ov 1-\al \z}\right)^n\left({\z-\al \ov \z+\al }\right)^m
z^{-i-1}\,\z^{-j-1}\,{\df z\,\df\z\ov z\z-1},\label{H1H2ij}\ee
where the contours are such that $|z\z|>1$.
Setting $i=j=0$ and making the substitution $\z\ra\z\inv$ gives
\[h=(H_1H_2)_{0,0}=\left({1\ov 2\pi \mi}\right)^2\intint
\left({1+\al z\ov 1-\al z}\right)^m\left({z-\al \ov z+\al }\right)^n
\left({\z+\al\ov \z-\al}\right)^n\left({1-\al\z \ov 1+\al\z }\right)^m
\,{\df z\,\df\z\ov z(z-\z)},\]
where now on the contours (both still described counterclockwise) $|z|>|\z|$.
If we ignore the $z-\z$ in the denominator we have two integrals to each of
which we apply steepest descent. The saddle points are $\pm \mi \beta$ for both.
The new $z$ contour is as before but the new $\z$ contour closes at $-\al$ and $\al\inv$.
We first deform the original $\z$ contour to this, always remaining inside the
original $z$ contour. Then we deform the $z$ contour to its steepest descent curve. In the
process we pass through the points of the $\z$ contour
from $-\mi\beta$ to $\mi\beta$ in the right half-plane. The $z$ residues at these points
are $1/\z$ and so the
deformations lead to the double integral over the steepest descent contours, which
is $O(n^{-1/2})$, plus $(2\pi i)\inv\int_{-\mi\beta}^{\mi\beta}d\z/\z=\hf$. This establishes
the second limit.

When $\om=1-\sqrt2$ we find that $\det S$ has the limit $(5-3\sqrt2)/8\ne0$ as
$t\to0,$\linebreak $h\to 1/2$.\footnote{Computations verify that when $\tau<\al^2$ the limit
of $h$ is 0 and $t=(-1)^m+o(1)$, and that the limit of $\det S$ is $(9-4\sqrt2)/2$
as $m\to\iy$ through even values and $-1/2$ as $m\to\iy$ through odd values. For
the limits when $\tau>\al^{-2}$ we replace $m$ by $n$. The proofs should be similar
to what we have already done, except that the saddle points will now be real.}
Hence the entries of $S\inv$ are all bounded. (The invertibility
of $S$ for large $n$ implies in turn the invertibility of $I+\cK$ and hence the validity of
the determinant identities we have been using.)

\subsection{Scaling}

If we write $P$ for $I-P_h$ then we see from (\ref{Kinvrep}) that $P(I+\cK)\inv P$ equals 
the identity operator $I$ minus
\[\twotwo{PH_1H_2P}{0}{0}{PH_2H_1P}-
\twotwo{0}{P(I-H_1H_2)L_x\La P}{P(I-H_2H_1)L_y\La P}{0}\]
\be +\sum_{i,j=1}^4 s_{ij}\,Pa_i\ot Pb_j'.\label{ab'sum}\ee

Eventually we will set $i=h+n^{1/3}x,\ j=h+n^{1/3}y$ where $h=cn+n^{1/3}s$ with $c$ to
be determined. The operators $PH_1$ and $H_1P$ will give rise to integrals like
\[\int\left({1+\al z\ov 1-\al z}\right)^m\left({z-\al \ov z+\al }\right)^n
z^{-cn-n^{1/3}x}\,\df z\]
(with a different $x$) and $PH_2$ and $H_2P$ will give rise to integrals like
\[\int\left({1+\al z\ov 1-\al z}\right)^n\left({z-\al \ov z+\al }\right)^m
z^{-cn-n^{1/3}x}\,\df z.\]
If we make the substitution $z\to z\inv$ in the latter we get an integral like
\[\int\left({z+\al\ov z-\al}\right)^n\left({1-\al z\ov 1+\al z}\right)^m
z^{cn+n^{1/3}x}\,\df z.\]
If we set
\[\psi(z)=\left({1+\al z\ov 1-\al z}\right)^m\left({z-\al \ov z+\al }\right)^n
z^{-cn},\]
our integrals become
\[\int \psi(z)\,z^{-n^{1/3}x}\,\df z \And \int\psi(z)\inv\,z^{n^{1/3}x}\,\df z.\]
If we think of the factors $\psi(z)^{\pm 1}$ as the
dominant ones and apply steepest descent, there will in general be two saddle
points for the two integrals, and the product of
the critical values will be exponentially small or large. If $c$ is chosen so the
critical points coincide\footnote{This $c$ will be the $c_1(\al,\tau)$ of the Theorem.}
then the product of the critical values will be 1 and the
product of the operators will have nontrivial scaling.

To determine $c$, let $\s(z)=n\inv\,\log\psi(z)$, so that
\[\s'(z)={2\al \tau\ov 1-\al^2z^2}+{2\al\ov z^2-\al^2}-{c\ov z}.\]
If we eliminate $c$ from $\s'(z)=\s''(z)=0$ we obtain
\be{\tau(1+\al^2z^2)\ov(1-\al^2z^2)^2}-{\al^2+z^2\ov(z^2-\al^2)^2}=0.\label{tau}\ee
The function on the left is strictly increasing from $-\iy$ to $+\iy$ on
the interval $(\al,\al\inv)$. It follows that there is a unique point $z_0$ in
this interval where the function vanishes. This will be our saddle point and we
set
\be c=2\al z_0\left({\tau\ov 1-\al^2z_0^2}+{1\ov z_0^2-\al^2}\right)>0.\label{c}\ee

{}From the behavior of $\s'(z)$ for large negative $z$ and near $-\al\inv,\ \al$ and 0
we see that $\s'$ has a zero in $(-\iy,-\al\inv)$ and a zero in $(-\al,0)$.
Since it has a double zero at $z_0$ this accounts for all four of its finite zeros.
Since $\s'(z)$ tends to $+\iy$ at the endpoints of $(\al,\al\inv)$ it follows that
it is positive everywhere there except at $z_0$, and this implies that
$\s'''(z_0)>0$. (Well, this only shows that $\s'''(z_0)\ge 0$.
We find in the last section an explicit expression for $\s'''(z_0)$
in terms of $z_0$, from which it is clear that it is positive.)

\begin{figure}
\bc
\resizebox{7cm}{7cm}{\includegraphics{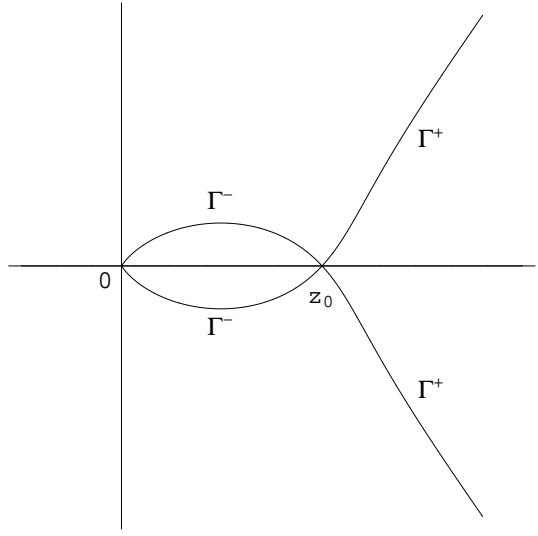}}
\caption{Steepest descent curves $\Ga^{\pm}$.}
\label{sdFig}
\ec
\end{figure}
The two steepest descent curves, which we call $\Ga^+$ for the first integral and $\Ga^-$
for the second, together form a single contour. The first emanates from $z_0$
at angles $\pm\pi/3$ with branches going to $\iy$ in two directions. The second emanates
from $z_0$ at angles $\pm 2\pi/3$ and closes at $z=0$.  See Figure~\ref{sdFig}.

It will be convenient to replace $\ell^2([h,\iy))$ by $\ell^2(\bZ_+)$, and to do that
we change our
meaning of the operator $P$. A $P$ appearing on the left
is to be interpreted as $\La^h$ and a $P$ appearing on the right is to be
interpreted as ${\La^*}^h$. So, for example, the $i^{\scriptstyle\mathrm{th}}$ component of $PT_1e$
is $(T_1e)_{h+i}$ and the $i,j$ entry of $PH_1H_2$ is $(H_1H_2)_{h+i,h+j}$. With
these reinterpretations of $P$ the operators in (\ref{ab'sum}) all act on $\ell^2(\bZ_+)$
and the determinant has not changed. (Recall that $h=cn+n^{1/3}s$.)

Let $D$ be the diagonal matrix with $i^{\scriptstyle\mathrm{th}}$  diagonal element equal to
$\psi(z_0)\inv z_0^{n^{1/3}s+i}$, and multiply (\ref{ab'sum}) by
$\left(\begin{array}{cc}D&0\\0&D\inv\end{array}\right)$
on the left and by
$\left(\begin{array}{cc}D\inv&0\\0&D\end{array}\right)$
on the right. This will not affect the determinant. (The reason is that all the
operators in (\ref{ab'sum}) have $i,j$ entry $O(r^{i+j})$ for fixed $n$, where $r$ can be 
any number
larger than $\al$, and $\al<z_0<\al\inv$.)
We shall show that after these multiplications the first operator scales
to the the direct sum of two Airy operators,\footnote{If
a matrix acting on $\ell^2(\bZ_+)$ has entries $M(i,j)$ and if the kernel
$n^{1/3}M([n^{1/3}x],[n^{1/3}y])$ acting on $L^2(0,\iy)$, which is unitarily equivalent to 
the
matrix operator, converges in trace norm to a limiting kernel
then we say that the matrix scales in trace
norm to the limiting kernel. The Fredholm determinant of $M(i,j)$ then converges
to the Fredholm determinant of the limiting kernel. This is what will happen here, with
the limiting kernel being the direct sum of two Airy kernels.} the second operator has 
trace norm $o(1)$ and the
vectors in the sum will all have norm $o(1)$.

\subsubsection{The first operator in (\ref{ab'sum})}

The upper left corner of the first operator becomes, aside from the identity operator,
$DPH_1H_2PD\inv$. We write this as
\[(DPH_1D_0)\;(D_0\inv H_2PD\inv),\]
where $D_0$ is the diagonal matrix with $i$th diagonal element equal to $z_0^{i}$\,,
and we scale each factor. We have
\[z_0^{n^{1/3}s+i+j} (PH_1)_{i,j}={z_0^{n^{1/3}s+i+j}\ov 2\pi \mi}\int
\left({1+\al z\ov 1-\al z}\right)^m\left({z-\al \ov z+\al }\right)^n
\,z^{-h-i-j-2}\,\df z\]
\be ={z_0^{n^{1/3}s+i+j}\ov 2\pi \mi}\int\psi(z)
\,z^{-n^{1/3}s-i-j}\,{\df z\ov z^2}.\label{H1integral}\ee

The main fact will be the following. Define
\[\psi(z,\ga)=\left({1+\al z\ov 1-\al z}\right)^m\left({z-\al \ov z+\al }\right)^n
z^{-\ga n}\]
and set
\[ I_n(x)={1\ov 2\pi \mi}\int
\psi(z,c+n^{-2/3}x)\,{\df z\ov z^2},\]
the contour being the unit circle. Then
\be\psi(z_0)\inv\,\,z_0^{n^{1/3}x}\,n^{1/3}\,|I_n(x)|\le e^{-\dl x}\label{estimate}\ee
valid for some $\dl$ and all $n$ if $x$ is bounded from below, and
\be\lim_{n\to\iy}\psi(z_0)\inv\,\,z_0^{n^{1/3}x}\,n^{1/3}\,I_n(x)=z_0^{-1}\,g\,
{\rm Ai}(gx)\label{limit}\ee
pointwise, with $g$ the constant given by (\ref{g}). (The limit in (\ref{limit}) will
be uniform for $x$ in a bounded set.)

We first show that (\ref{estimate}) holds if $n^{1/3}x>\eta n$  for some $\eta>0$. For this
we set $\ga=c+n^{-2/3}x$ so that $\ga-c>\eta$, write our main integrand as
$\psi(z,\ga)$ and do steepest descent. With
\be \s(z,\ga)=n\inv\log\psi(z,\ga)=\s(z)-(\ga-c)\log z,\label{szga}\ee
our saddle points $z_{\ga}^{\pm}$ (there will be two of them when $\ga>c$) are solutions of
$\s'(z_\ga^{\pm},\ga)=0.$
Differentiating this with respect to $\ga$ gives
\[\s''(z_\ga^{\pm},\ga){dz_\ga^{\pm}\ov d\ga}-{1\ov z_\ga^{\pm}}=0.\]
Since $\s''(z_\ga^{\pm},\ga)\ne0$ (since $c$ is the only value of $\ga$ for which there
is a double saddle point) we find that if $z_\ga^{\pm}$ are chosen so that
$z_\ga^+>z_\ga^-$ for $\ga$ near $c$ then $z_\ga^+$
increases always as $\ga$ increases and $z_\ga^-$ decreases, and $\s_1''(z_\ga^+,\ga)>0$
and $\s_1''(z_\ga^-,\ga)<0$. In particular $z_\ga^+$ is the saddle point we take for our steepest
descent, and $z_\ga^+>z_0$ when $\ga>c$ since $z_c^+=z_0$. For the critical value
we have to see how $\s(z_\ga^+,\ga)$  behaves as a function
of $\ga$. From (\ref{szga}) and the fact $\s'(z_\ga^\pm,\ga)=0$ we obtain
\[{d\s(z_\ga^\pm,\ga)\ov d\ga}=-\log z_\ga^\pm,\]
and so
\be{d\ov d\ga}[\s(z_\ga^\pm,\ga)+\ga\log z_0]=\log {z_0\ov z_\ga^\pm}.\label{dsdgpm}\ee
Since $z_\ga^+$ is an increasing function of $\ga$ this shows that
$\s(z_\ga^+,\ga)+\ga\log z_0$ is a decreasing
concave function of $\ga$. It follows that for some $\dl>0$
\[\s(z_\ga^+,\ga)+\ga\log z_0< \s(z_0)+c\log z_0-\dl\,(\ga-c)\]
for $\ga\ge c+\eta$. Thus for these $\ga$
\[\psi(z_0)\inv\,z_0^{n(\ga-c)}\,\psi(z_\ga^+,\ga)\le \me^{-\dl n(\ga-c)}.\]
Since $|\psi(z,\ga)|$ achieves its maximum on the steepest descent contour at $z_\ga$,
and since the contour is bounded away from zero (this follows from the fact that $z_\ga
\to\al$ as $\ga\to\iy$), we see that in this case $I_n(x)$ is at most a constant times
$\psi(z_\ga^+,\ga)$, where $\ga=c+n^{-2/3}x$. Hence
\[\psi(z_0)\inv\,\,z_0^{n^{1/3}x}\,|I_n(x)|=O(\me^{-\dl n^{1/3}x}),\]
when $n^{1/3}x>\eta n$. This is an even better estimate than
(\ref{estimate}).

Since we have shown that (\ref{estimate}) holds if $n^{1/3}x>\eta n$, where $\eta$
can be as small as we please, we may assume $n^{1/3}x=o(n)$ when $x>0$.

Write
\[I_n(x)={1\ov 2\pi \mi}\int\psi(z)\,z^{-n^{1/3}x}\,{\df z\ov z^2},\]
and use the steepest descent curve $\Ga^+$. It emanates from $z_0$ at
angles $\pm \pi/3$. Clearly it is bounded away from 0. Choose $\ve$ small
and let
\[\Ga^{(1)} =\{z\in\Ga^+:|z-z_0|>\ve\} \And \Ga^{(2)} =\{z\in\Ga^+:|z-z_0|<\ve\}\]
with corresponding $I_n^{(1)}(x)$ and $I_n^{(2)}(x)$. We shall show that both of these
satisfy the uniform estimate (\ref{estimate}) and
\[\lim_{n\to\iy}\psi(z_0)\inv\,n^{1/3}\,I_n^{(1)}(x)=0\And
\lim_{n\to\iy}\psi(z_0)\inv\,z_0^{n^{1/3}x}\,n^{1/3}\,I_n^{(2)}(x)=
z_0^{-1}\,g\,{\rm Ai}(gx).\]

Consider $I_n^{(1)}(x)$ first. Since $\Re\,\s(z)$ is strictly decreasing as we move
away from $z_0$ on $\Ga$ we know that
$|\psi(z)|< \psi(z_0)\,\me^{-\dl n}$ on $\Ga^{(1)}$ for some $\dl>0$. The assertions follow
from this since we are in the case $n^{1/3}x=o(n)$. (This also holds also for $x<0$ since
$x$ is bounded below.)

Making the variable change $z\to z_0\,(1+\x)$ we can write, since $z=z_0\,e^{\x(1+O(\x))}$
near $z=z_0$,
\[ I_n^{(2)}(x)=\psi(z_0)\,z_0^{-n^{1/3}x}{1\ov2\pi \mi}\int_{|\x|<\ve/z_0}
\me^{nbz_0^3\x^3(1+O(\x))-n^{1/3}x\x(1+O(\x))}\,(z_0^{-1}+O(\x))\,\df\x,\]
where $b=\s'''(z_0)/6$. The path of integration here is the portion of the contour
$\Ga^{(2)}$ satisfying the indicated
inequality. It consists of two little arcs emanating from $\x=0$ tangent to the line
segments making angles $\pm \pi/3$ with the positive axis. If $\ve$ is small enough
and we replace the integral by
the line segments themselves we introduce an error of then form $O(e^{-\dl n})$ with
a different $\dl$, since $n^{1/3}x=o(n)$.
With the variable change $\x\to n^{-1/3}\x$ we obtain
\[\psi(z_0)\inv\,z_0^{n^{1/3}x}\,n^{1/3}\,I_n^{(2)}(x)\]
\[={1\ov2\pi \mi}\int_{|\x|<n^{1/3}\ve}
\me^{bz_0^3\x^3(1+O(n^{-1/3}\x))-x\x(1+O(n^{-1/3}\x))}(z_0^{-1}+O(n^{-1/3}\x))\,\df\x+O(e^{-\dl n}),\]
where now the integration is taken over line segments of length of
the order $n^{1/3}$.

On the path of integration we have $\Re(\x^3)\le-\dl|\x|^3$ and $\Re(\x)\ge\dl|\x|$ for
some $\dl>0$ and so for the above we have an estimate of the form
\[\int_0^\iy \me^{-\dl (t^3-xt)}\df t+O(e^{-\dl n}).\]
Since $n\gg x$ this gives the required uniform bound. The limit
of the integral,
with its factor $1/2\pi \mi$, equals $z_0^{-1}\,g\,{\rm Ai}(gx)$, where
\be g=(3bz_0^3)^{-1/3}=z_0\inv\left({2\ov \s'''(z_0)}\right)^{1/3}.\label{g}\ee
That the limit is as
stated follows by taking the limit under the integral sign, which
is justified by dominated convergence.

To obtain the scaling of the matrix $DPH_1D_0$ we need only observe that by
(\ref{H1integral})
\[n^{1/3}\,(DPH_1D_0)_{[n^{1/3}x],\,[n^{1/3}y]}=
\psi(z_0)\inv\,\,z_0^{n^{1/3}s+[n^{1/3}x]+[n^{1/3}x]}\,n^{1/3}\,I_n(s+n^{-1/3}([n^{1/3}x]+[n^{1/3}y])).\]
It follows from (\ref{estimate}) and (\ref{limit}) that this kernel on $(0,\iy)$
converges in Hilbert-Schmidt norm to $z_0^{-1}\,g\,{\rm Ai}(g(s+x+y))$.

To scale $D_0\inv H_2PD\inv$ we write
\[z_0^{-n^{1/3}s-i-j} (PH_2)_{i,j}={z_0^{-n^{1/3}s-i-j}\ov 2\pi \mi}\int
\left({1+\al z\ov 1-\al z}\right)^n\left({z-\al \ov z+\al }\right)^m
\,z^{-cn-n^{1/3}s-i-j-2}\,\df z.\]
If we make the substitution $z\to z\inv$ this becomes
\[{z_0^{-n^{1/3}s-i-j}\ov 2\pi \mi}\int\psi(z)\inv\,z^{n^{1/3}s+i+j}\,\df z.\]
This is completely analogous to (\ref{H1integral}). To use the analogous argument we
mention only that we use (\ref{dsdgpm}) with the minus signs to see that both
$-\s(z_\ga^-,\ga)-\ga\log z_0$ and its derivative are decreasing functions of $\ga$.
The steepest descent curve now is $\Ga^-$. We need not go through the details again.
We find that
\[n^{1/3}\,(D_0\inv PH_1D\inv)_{[n^{1/3}x],\,[n^{1/3}y]} \]
converges in Hilbert-Schmidt norm to $z_0\,g\,{\rm Ai}(g(s+x+y))$.
Hence
\[n^{1/3}(DPH_1H_2PD\inv)_{[n^{1/3}x],\,[n^{1/3}y]}\]
converges in trace norm to $g\,K_{{\rm Airy}}(g(s+x),\,g(x+y))$ on $L^2(0,\iy)$.

By taking transposes we see that $D\inv PH_2H_1PD\inv$ has the same scaling limit, which
takes care of the lower right corner of the first operator in (\ref{ab'sum}).

\subsubsection{The second operator in (\ref{ab'sum})}

Next, we have to look at
\[D\, P(I-H_1H_2)L_x\La P\,D\ \ \textrm{and}\ \
D\inv\,P(I-H_2H_1)L_y\La P\,D\inv.\]
We find, using (\ref{HHL}),
\[z_0^{2n^{1/3}s+i+j}((I-H_1H_2)L_x\La)_{h+i,h+j}\]\[=\left({1\ov 2\pi \mi}\right)^2\intint
\left({1+\al z\ov 1-\al z}\right)^m\left({z-\al \ov z+\al }\right)^n
\left({1+\al \z\ov 1-\al \z}\right)^n\left({\z-\al \ov \z+\al }\right)^m\]
\[\times z_0^{2n^{1/3}s+i+j}{z^{-h-i-1}\,\z^{h+j-1}\ov z\z-1}\,\df z\,\df\z\,(-1)^{h+j}.\]
After the substitution $\z\to-\z\inv$ this becomes
\[\left({1\ov 2\pi \mi}\right)^2\intint
\psi(z)\,\psi(\z)\,\left({z\ov z_0}\right)^{-n^{1/3}s-i}
\left({\z\ov z_0}\right)^{-n^{1/3}s-j}
{\df z\,\df\z\ov z(z+\z)}.\]
The integrals here are initially taken over circles close to the unit circle,
with $|z|>|\z|$.
We first deform the $\z$ contour to $\Ga^+$, while always
having $z+\z$ nonzero. Then if we deform the $z$ contour we pass through a pole
at $\z=-z$ for every $z\in\Ga^+$. The residue at the pole equals a constant times
$z^{-2h-i-j-2}$, and integrating this over $\Ga^+$ gives zero. So
both integrals may be taken over $\Ga^+$.
Since now the denominator does not vanish at $z=\z=z_0$ the same sort of argument
we already gave shows that
this operator equals a constant times $\psi(z_0)^2\,n^{-1/3}$
times an operator which scales to the trace class operator
$g\,{\rm Ai}(g(s+x))\ot g\,{\rm Ai}(g(s+y))$. In particular its trace norm is
$O(\psi(z_0)^2\,n^{-1/3})$. This shows that $D\, P(I-H_1H_2)L_x\La P\,D$ has
trace norm $O(n^{-1/3})$ and a similar argument applies to
$D\inv\,P(I-H_2H_1)L_y\La P\,D\inv$.

\subsubsection{The last operator in (\ref{ab'sum})}

Finally we consider the vectors $Pa_i$ and $Pb_j'$ and look at their constituents
$PT_1e,\ PT_2e$,   $H_1H_2e$ and $H_2H_1e$. We shall show that
\[\|DPT_1e\|=O(n^{-1/6}) \And |D\inv PT_2e\|=O(n^{-1/6}),\]
\[\|DH_1H_2e\|=O(n^{-2/3}) \And \|D\inv H_2H_1e\|=O(n^{-2/3}).\]

For the first, we have
\[(T_1e)_{h+i}={1\ov2\pi \mi}\int\psi(z)\,z^{-n^{1/3}s-i-1}\,\df z=I_n(s+n^{-1/3}(i-1)),\]
and from (\ref{estimate}) and (\ref{limit}) we deduce now that the function
\[n^{1/3}(DPT_1e)_{[n^{1/3}x]}\]
converges in $L^2(0,\,\iy)$. In particular its norm is $O(1)$. But then
\[\|\{(DPT_1e)\}_i\|_{\ell^2}=n^{1/6}\,\|(DPT_1e)_{[n^{1/3}x]}\|_{L^2}=O(n^{-1/6}).\]
Similarly $\|D\inv PT_2e\|=O(n^{-1/6})$.

For $DH_1H_2e$, we have from (\ref{H1H2ij})
\[z_0^{n^{1/3}s+i}\,(H_1H_2e)_{h+i}\]
\[=\left({1\ov 2\pi \mi}\right)^2\intint
\left({1+\al z\ov 1-\al z}\right)^m\left({z-\al \ov z+\al }\right)^n
\left({1+\al \z\ov 1-\al \z}\right)^n\left({\z-\al \ov \z+\al }\right)^m
\,z^{-h-i-1}\,z_0^{n^{1/3}s+i}\,\z^{-1}\,{\df z\,\df\z\ov z\z-1},\]
and with the substitution $\z\to\z\inv$ this becomes
\be \left({1\ov 2\pi \mi}\right)^2\intint
\psi(z)\left({\z+\al\ov \z-\al}\right)^n\left({1-\al\z \ov 1+\al \z}\right)^m
z^{-n^{1/3}s-i-1}\,z_0^{n^{1/3}s+i}
\,{\df z\,\df\z\ov z-\z}.\label{DHH}\ee
As before the integrals here are initially taken over circles close to the unit circle,
with $|z|>|\z|$. Now we want to deform the $z$ contour
to $\Ga^+$ and the $\z$ contour to its steepest descent contour $C$, a curve passing
through the saddle points $\pm i\beta\inv$ and closing at $-\al$ and $\al\inv$.

To do this we show first that, except for $z_0$, all points of $\Ga^+$ satisfy $|z|>z_0$.
This will follow if we can show that on the circle $z=z_0\,\me^{\mi\th}$ the absolute minimum of
\be\log\left|\left({1+\al z\ov1-\al z}\right)^\tau\left({z-\al\ov z+\al}\right)\right|
\label{logabs}\ee
occurs at $\th=0$. (For then $|\psi(z)|$ would be larger than $\psi(z_0)$ everywhere
on the circle except for $z=z_0$, so no other point on the circle could be on $\Ga^+$.
Locally $\Ga^+$ is outside the circle and so it would have to be everywhere outside.)
Using
\[{d\ov d\th}=\mi z\,{d\ov dz}\]
we find that the derivative with respect to $\th$ of (\ref{logabs}) equals
$-2\al$ times the imaginary part of
\[{\tau\ov z\inv-\al^2z}+{1\ov z-\al^2z\inv}={(\tau-\al^2)z+(1-\tau\al^2)z\inv
\ov (z\inv-\al^2z)\,(z-\al^2z\inv)}.\]
This vanishes exactly when the imaginary part of
\[((\tau-\al^2)z+(1-\al^2)z\inv)\,(\bar{z}\inv-\al^2\bar{z})\,(\bar{z}-\al^2\bar{z}\inv)\]
does. This is a trigonometric polynomial in $\th$ of degree 3. It is an odd function of
$\th$ and so is of the form $\sin\th$ times a polynomial of degree two in $\cos\th$. Since
it has at least a double zero at $\th=0$ (by the choice of $c$ and $z_0$) the polynomial
must have a factor $\cos\th-1$. Since it is
an odd function of $z$ it must also have a double zero at $\th=\pi$, so there must also
be a factor $\cos\th+1$. Thus it must be equal to a constant times $\sin\th\,(\cos^2\th-1)$.
In particular there can be no other zeros. Thus (\ref{logabs}), which we know has a local minimum
on the circle at $z=z_0$, must have its absolute minimum there (and its absolute
maximum at $z=-z_0$). Thus, as claimed, all points of $\Ga^+$ except for $z_0$ satisfy $|z|>z_0$

In particular, all points of $\Ga^+$ satisfy $|z|\ge z_0$. Since $z_0>\al$ we can
first take the integrals in (\ref{DHH}) over the circles $|z|=z_0$ and $|\z|=\al+\ve$
with $\ve$ small and positive. Then we can deform the $z$ contour to $\Ga^+$ without
crossing the circle $|\z|=\al+\ve$. Next we want to deform the $\z$ contour to $C$.
This curve closes on the right at $\al\inv$ and so, since $z_0<\al\inv$, it intersects
$\Ga^+$ at two points $z'$ and $z''$, say. (In principle there could
be finitely many other points; the following argument could be easily modified in this case.)
Hence upon deforming the $\z$ contour to $C$ we
pass through a pole for those $z$ on the arc of
$\Ga^+$ passing through $z=z_0$ with end-points $z'$ and $z''$. For each $z$ on this
arc the residue at $\z=z$ equals $z^{-h-i-1}\,z_0^{n^{1/3}s+i}$, and then integrating with
respect to $z$ gives
\[z_0^{n^{1/3}s}\,(h+i)\inv \left[{z_0^i\ov {z''\,}^{h+i}}-{z_0^i\ov {z'\,}^{h+i}}\right].\]
We claim that $\psi(z_0)\inv$ (which is the factor contained in $D$) times
this vector is exponentially small, i.e., $O(\me^{-\dl n})$ for some $\dl>0$.

Because all points of $\Ga^+$ except for $z_0$ satisfy $|z|>z_0$ the vectors
$\{(z_0/z'')^i\}$ and $\{(z_0/z')^i\}$ belong to $\ell^2(\bZ_+)$.
So we need only show that $\psi(z_0)\,\inv|z'|^{-h}$ is
exponentially small (and the same for $z''$).
In fact, on the part of $C$ in the right half-plane
\[\left|\left({z+\al\ov z-\al}\right)\left({1-\al z\ov 1+\al z}\right)^\tau\right|\]
is at most 1 and is strictly less than 1 outside a neighborhood of the critical points
$\pm i\beta$. This shows that
$|z'|^{-c}\le (1-\dl)\,|\psi(z')|^{1/n}$
for some $\dl>0$. Since $|\psi(z)|<|\psi(z_0)|$ on $\Ga^+$ this shows that
$|z'|^{-c}\le (1-\dl)\,|\psi(z_0)|^{1/n}$ and it follows that
$\psi(z_0)\inv\,|z'|^{-h}$ is exponentially small.

We can now say that, with error $\psi(z_0)^2$ times an exponentially small quantity,
the square of the norm of the vector (\ref{DHH}) equals a quadruple integral in which
every term in the integrand
except $z^{-h-i-1}\,z_0^i$ has an analogous term with variables $z',\ \z'$, and the
$z^{-h-i-1}\,z_0^i$
term becomes $(zz')^{-h}/(zz'-z_0^2).$  For the $z$ and $z'$ integrals we integrate over
$\Ga^+$ and for the $\z$ and $\z'$ integrals we integrate over $C$. (We use again here
the fact that $|z|,\ |z'|>z_0$, so we can sum under the integral signs.) The
$\z,\ \z'$ integrals contribute $O(n\inv)$
while the $z,\ z'$ integrals contribute $O(\psi(z_0)^2\,n^{-1/3})$. Thus
$\|DPH_1H_2e\|^2=O(n^{-4/3})$.\footnote{This was under the basic assumption $\al^2<\tau<
\al^{-2}$. Otherwise the $\z$ and $\z'$ integrals are only $O(1)$, with the result that
$\|DH_1H_2e\|$ is only $O(n^{-1/6})$ and the same will hold for $\|D\inv H_2H_1e\|$. 
These are still good enough since $o(1)$ is all that is needed.} 

{}For $H_2H_1e$ we interchange $m$ and $n$ and make the variable change $z\to z\inv$ but not
the variable change $\z\to\z\inv$. Thus
\[z_0^{-n^{1/3}s-i}\,(H_2H_1e)_{h+i}\]
\[=\left({1\ov 2\pi \mi}\right)^2\intint
\psi(z)\inv
\left({1+\al \z\ov 1-\al \z}\right)^n\left({\z-\al \ov \z+\al }\right)^m
\,z^{n^{1/3}s+i}\,z_0^{-n^{1/3}s-i}\,{\df z\,\df\z\ov \z-z}.\]
Here originally we must have $|\z|>|z|$ on the contours and we want to deform them so
that the $\z$ contour becomes $\Ga^-$ and the $\z$ contour becomes what we again call $C$.
Now in the deformation, we pass through a pole in the $\z$ integration for those $z$
on an arc of $\Ga^-$ passing through $z=z_0$ with end-points which we again
call $z'$ and $z''$. The residue equals $z^{h+i}\,z_0^{-n^{1/3}s-i}$
and integration with respect to $z$ gives
\[z_0^{-n^{1/3}s}\,(h+i+1)\inv \left[{{z''\,}^{h+i+1}\ov z_0^i}-
{{z'\,}^{h+i+1}\ov z_0^i}\right].\]
This is completely analogous to what went before. Now $|z'|,\ |z''|<z_0$ and
$\psi(z_0)\,|z'|^{h}$ is exponentially small. We continue
as with $H_1H_2e$ and find that $\|D\inv PH_2H_1e\|^2=O(n^{-4/3})$.

If we go back to the forms of the $a_i$ and $b_j'$ described
earlier we see that the vectors which arise
after multipliying by the diagonal matrices $D$ and $D\inv$ are exactly the four whose
norms we just estimated.

\subsection{Recapitulation}

We have shown that the matrix (\ref{ab'sum}) acting on
$\ell^2(\bZ_+)$ scales in trace norm to the kernel
\[\twotwo{g\,K_{{\rm Airy}}(g(s+x),\,g(x+y))}{0}{0}{g\,K_{{\rm Airy}}(g(s+x),\,g(x+y))}\]
acting on $L^2(0,\iy)$. It follows that its Fredholm determinant converges to $F(gs)^2$. In
view of (\ref{probrep}) and   (\ref{Kdetrep})
this establishes that for fixed $s$
\[\lim_{n\to\iy}\pr_{\sigma}\,(L\le c\,n+n^{1/3}s)=F_2(gs),\]
where $c$ is determined by (\ref{tau}) and (\ref{c}) and $g$ by (\ref{g}).
This gives the statement of the Theorem, where the constants $c_1(\al,\tau)$ and
$c_2(\al,\tau)$ of the Introduction
are, respectively, $c$ and $g^{-1}$.

\subsection{Computation of $\mathbf{\s'''(z_0)}$}

Think of $c$ and $z_0$ as functions of $\tau$, which they are. We have
\[\s''(z)=4\al z\left[{\al^2\tau\ov (1-\al^2z^2)^2}-{1\ov(z^2-\al^2)^2}\right]+
{c\ov z^2}.\]
Differentiating the identity $0=\s''(z_0)$ with respect to $\tau$ and using the above
give
\be 0=\s'''(z_0)\,z_0'+{4\al^3\,z_0\ov (1-\al^2z_0^2)^2}+{c'\ov z_0^2},\label{s'''}\ee
where $z_0'$ and $c'$ denote $dz_0/d\tau$ and $dc/d\tau$, respectively.

{}From (\ref{tau}) and (\ref{c}) we
find that $c$ is given in terms of $z_0$ by the relation
\[c={4\,\al\,(1-\al^4)\,z_0^3\ov(1+\al^2z_0^2)\,(z_0^2-\al^2)^2}.\]
We compute that
\be{dc\ov dz_0}=-4{\al\,(1-\al^4)\,z_0^2\,(z_0^2+3\al^2+3\al^2z_0^4+\al^4z_0^2)
\ov(1+\al^2z_0^2)^2\,(z_0^2-\al^2)^3}.\label{dcdz}\ee

{}From (\ref{tau}) $\tau$ is given in terms of $z_0$ by
\[\tau={(1-\al^2z_0^2)^2\,(\al^2+z_0^2)\ov(1+\al^2z_0^2)\,(z_0^2-\al^2)^2}.\]
We compute that
\be{d\tau\ov dz_0}=-{2\,(1-\al^4)\,z_0\,(1-\al^2z_0^2)\,(\al^4z_0^2+3\al^2z_0^4+3\al^2+z_0^2)
\ov(z_0^2-\al^2)^3\,(1+\al^2z_0^2)^2}.\label{dtdz}\ee

We first solve (\ref{s'''}) for $\s'''(z_0)$ in terms of $z_0'$ and $c'$ (and $z_0$).
Then we use
(\ref{dcdz}) and (\ref{dtdz}) and the
relations $z_0'=(d\tau/dz_0)\inv$ and $c'=(dc/dz_0)\,z_0'$.
We find that
\[\s'''(z_0)=4{\al\,(1-\al^4)\,[(1+\al^4)\,z_0^2+3\al^2\,(1+z_0^4)]\ov (1-\al^4z_0^4)
\,(z_0^2-\al^2)^3}.\]
This is positive since $0<\al<1$ and $\al<z_0<\al\inv$.\sp

\section{Poisson Limit of the Shifted Schur Measure\label{refereeSec}}
\setcounter{equation}{0}

For Schur measure  there are two interesting limiting cases:   the exponential limit
and the Poisson limit.  The exponential
limit of shifted Schur measure
 is supported on standard shifted tableaux; and hence, it is
expressible in terms of $f^{\la}_s$ (recall (\ref{shiftedSYT})).
The Poisson limit of
shifted Schur measure yields a natural interpretation of the maximizing rule
as a symmetry condition of the process.\footnote{The following remarks
are due to the referee.}  Namely, if one sets
\[ m=n, \quad \alpha={t\ov n} \]
and takes $n\ra\iy$, the percolation-type model described in \S\ref{RSKsec}
becomes the following:  Consider two Poisson processes of rate $t^2$ both in
the (same) square $[0,1]\times [0,1]$.  Hence one can imagine two types of points,
marked and unmarked, in the square.  Now with probability one no two points,
whether marked or unmarked, have the same $x$ or $y$ coordinates; and hence,
the strictly increasing conditions (2) and (3) of\ \S\ref{RSKsec} for the rule
of the maximizing path are not necessary.  Therefore, we do not need to distinguish
the unmarked and marked points.  Since the union of two Poisson processes
has the rate $2t^2$, the reulting process is as follows:  In the square
$[0,1]\times [0,1]$ select Poisson points of rate $2t^2$.  Then take the longest
path starting from the lower right corner $(1,0)$ that follows an up/left
path, turns in direction once and only once, and follows an up/right path ending
at the upper right corner $(1,1)$.

It is clear that this length is also equal to the following symmetric version.
Take a realization of the Poisson process.  Take the mirror image of the points
about the left side.  Adjoin the mirror image on the left and the original points
on the right.  Hence the resulting rectangle has sides of lengths $2$ and $1$,
and there are twice as many points of the orginal configuration which are
symmetric about the center vertical line.  The (usual) longest up/right path
from the left lower corner to the right upper corner is precisely the longest
maximizing path from the lower right corner to the upper right corner in the
above description.

The (formal) limit $m=n\ra\iy$ with $\alpha=t/n$ ($t$ fixed) in the main theorem of \S 1
is
\begin{equation}
 \lim\pr\left({L-4t\ov (2t)^{1/3}}<s\right)=F_2(s). \label{symmetryFluc}\end{equation}
The consequence is that the vertically symmetric Poisson process has
the same fluctuation as the usual Poisson process with no symmetry condition.
Also the scaling in the above result is consistent with this intuition.
There are Poisson points of rate $4t^2$ (double of $2t^2$) in the square
of sides 2 and 1.  For such the case the limit (\ref{symmetryFluc}) is also
valid for the case of no symmetry condition.

\begin{center} {\bf Acknowledgments} \end{center}

This work was supported by the National Science Foundation through grants
DMS-9802122 and DMS-9732687. The authors thank Richard Stanley and Sergey
Fomin for useful early discussions concerning the RSK correspondence and
the referee for the remarks of \S\ref{refereeSec}. Finally, we
acknowledge our appreciation to the administration of the Mathematisches
Forschungsinstitut Oberwolfach for their hospitality during the authors' visit under
their Research in Pairs program.\sp

\end{document}